\documentclass[reqno,12pt]{amsart}

\usepackage{amsthm,amsmath,amssymb}
\usepackage{mathrsfs,amsfonts,dsfont,functan,extarrows,mathtools}
\usepackage[colorlinks,linkcolor=black,citecolor=red]{hyperref}

\usepackage{marginnote}
\usepackage{xcolor}

\makeatletter

\newcommand{\Rmnum}[1]{\expandafter\@slowromancap\romannumeral #1@}
\makeatother

\newtheorem{theorem}{Theorem}[section]
\newtheorem{proposition}[theorem]{Proposition}

\newtheorem{remark}[theorem]{Remark}
\newtheorem{lemma}[theorem]{Lemma}

\numberwithin{equation}{section}
\allowdisplaybreaks

\newcounter{wronumber}

\begin{document}

\title[]{K\"ahlerness of compact Hermitian surfaces under Semi-definite Strominger-Bismut-Ricci Curvatures}

\author[]{Liangdi Zhang}
\address[Liangdi Zhang]{\newline Mathematical Science Research Center, Chongqing University of Technology, Chongqing 400054, China}
\email{ldzhang91@163.com}



\begin{abstract}
We prove several K\"ahlerness criteria for compact Hermitian surfaces under semi-definiteness assumptions on natural Ricci curvatures of the Strominger-Bismut connection. The key tools for proving these results are explicit identities relating these Ricci curvatures to the torsion of the Strominger-Bismut connection, together with corresponding Chern number identities on compact Hermitian surfaces. The results may be viewed as Strominger-Bismut analogues and reformulations of Yang's K\"ahlerness criteria for compact complex surfaces.

\vspace*{5pt}

\noindent{\it Keywords}: Compact Hermitian surface; Strominger-Bismut-Ricci curvature; K\"ahler

\noindent{\it 2020 Mathematics Subject Classification}: 53C55
\end{abstract}

\maketitle

\tableofcontents

\section{Introduction}
\label{sec1}
Let $(M,\omega)$ be a compact Hermitian surface with $\omega=\sqrt{-1}h_{i\bar{j}}dz^i\wedge d\bar{z}^j$. The Strominger-Bismut connection (also known as Strominger connection or Bismut connection) first appeared in theoretical physics: Strominger \cite{Str86} introduced it in the study of heterotic string compactifications with torsion, where the torsion 3-form corresponds to the flux field strength in supersymmetric backgrounds. Independently, Bismut  \cite{Bis89} rediscovered the same connection in complex differential geometry, proving a local index theorem on non-K\"ahler manifolds by exploiting its favorable analytic properties. For a comprehensive account of this topic, we refer to  \cite{AV22}, \cite{AOUV22},  \cite{FOU15}, \cite{FTV22}, \cite{FU13}, \cite{LY12}, \cite{WY25}, \cite{WYZ20}, \cite{YZ18a}, \cite{YZ18b}, \cite{Yang2507}, \cite{Yang2508}, \cite{YZZ23}, \cite{Ye25}, \cite{ZZ19}, \cite{ZZ23}, \cite{ZZ24}, \cite{ZZ25} and the references therein.

The Enriques-Kodaira classification theorem (see \cite[Chapter \Rmnum{6}]{book}, \cite{Enri14,Enri49,Kod64,Kod66,Kod68a,Kod68b}) classifies nonsingular minimal compact complex surfaces into several classes according to their Kodaira dimension and other birational invariants, each parametrized by a moduli space. These classes fall into two broad types: K\"ahler surfaces, which include complex tori, $K3$ surfaces, and surfaces of general type, among others, and non-K\"ahler surfaces, which occur primarily in Class \Rmnum{7}. Belgun's refinement (\cite{Bel00}) shows that a compact non-K\"ahler admits a locally conformally K\"ahler metric with parallel Lee form if and only if it belongs to the corresponding class described in Belgun's classification. For the classification of non-K\"ahler surfaces, particularly those of Class \Rmnum{7}, via geometric flows, see for example \cite{Bo16}, \cite{GJS23}, \cite{S20}, \cite{ST10}, \cite{ST13}, \cite{Ye24} for approaches based on the pluriclosed flow, and for example \cite{Ed21}, \cite{TW13}, \cite{TW15}, \cite{TWY15} for those based on the Chern-Ricci flow.

On K\"ahler surfaces, curvature notions from Chern connection $^C\nabla$, induced Levi-Civita connection ${^{iLC}}\nabla$ (see \cite{LY17}) and Strominger-Bismut connection ${^{SB}}\nabla$ largely coincide, whereas on non-K\"ahler Hermitian surfaces the presence of torsion leads to diverse curvature behaviors, making them a natural testing ground for K\"ahlerness theorems under sign conditions on Hermitian Ricci curvatures.

It is well-known that a compact Hermitian manifold with positive (the first) Chern-Ricci curvature must be K\"ahler. In 2025, Yang \cite{Yang2508} established an explicit relation between the complexification of the real Ricci curvature of the complexified Levi-Civita connection $^{LC}\nabla$ and the torsion of Hermitian metrics. As an application, a compact Riemannian 4-manifold is a K\"ahler surface if it admits a compatible complex structure with vanishing $(2,0)$-component of the complexified Riemannian Ricci curvature and the $(1,1)$-component satisfies that $\mathfrak{R}ic^{(1,1)}+\frac{\sqrt{-1}}{4}\bar{\partial}^*\omega\wedge\partial^*\omega\leq0$, which in the Gauduchon case reduces to $\mathfrak{R}ic^{(1,1)}\leq0$. For systematic relations among the Ricci curvatures associated with Gauduchon connections on Hermitian manifolds, we refer to Broder-Stanfield \cite{BS} and Wang-Yang \cite{WY25}.

Yang \cite{Yang2508} also established Chern number identities on compact complex surfaces and show that a compact Riemannian four-manifold with constant Riemannian scalar curvature is K\"ahler if it admits a compatible complex structure such that the complexified Ricci curvature is a non-positive $(1,1)$-form. Motivated by Yang's work, we establish Strominger-Bismut counterparts of these K\"ahlerness criteria. The novelty lies in deriving explicit comparison identities between the complexified real Ricci curvature of the Strominger-Bismut connection, its natural Ricci contractions, and torsion terms, and then applying these identities to compact Hermitian surfaces. We emphasize that the main results below should not be interpreted as a direct strengthening of Yang's theorems. Rather, they are Bismut-type analogues obtained from curvature identities specific to the Strominger-Bismut connection.

Let $g$ be the background Riemannian metric and $J$ be the complex structure satisfying
\begin{equation}
g(X,Y)=g(JX,JY),\ \ \omega(X,Y)=g(JX,Y)
\end{equation}
for any $X$, $Y\in \Gamma(M,T_{\mathbb{R}}M)$, and
\begin{equation}
g(W,Z)=h(W,Z)
\end{equation}
for any $W$, $Z\in \Gamma(M,T_{\mathbb{C}}M)$ with $T_{\mathbb{C}}M=T_{\mathbb{R}}M\otimes\mathbb{C}=T^{1,0}M\oplus T^{0,1}M$.

We first fix our curvature conventions. The real curvature tensor of the Strominger-Bismut connection $^{SB}\nabla$ on the underlying Riemannian 4-manifold $(M,g,J)$ is defined by
\[R^{SB,\mathbb{R}}(X,Y,Z,W)=g({^{SB}}\nabla_X{^{SB}}\nabla_YZ-{^{SB}}\nabla_Y{^{SB}}\nabla_XZ-{^{SB}}\nabla_{[X,Y]}Z,W)\]
for any $X$, $Y$, $Z$, $W\in \Gamma(M,T_{\mathbb{R}}M)$.

The real Ricci curvature of $^{SB}\nabla$ on $(M,g,J)$ is defined by
\begin{equation}\label{1.0}
\mathcal{R}ic^{SB,\mathbb{R}}(X,Y)=g^{il}R^{SB,\mathbb{R}}\big(\frac{\partial}{\partial x^i},X,Y,\frac{\partial}{\partial x^l}\big)
\end{equation}
for any $X$, $Y\in \Gamma(M,T_{\mathbb{R}}M)$. In particular,
\[\mathcal{R}^{SB,\mathbb{R}}_{ij}=\mathcal{R}ic^{SB,\mathbb{R}}\big(\frac{\partial}{\partial x^i},\frac{\partial}{\partial x^j}\big)=g^{kl}R^{SB,\mathbb{R}}_{kijl}.\]

We denote by $\mathcal{R}ic^{SB,\mathbb C}$ the complex bilinear extension of $\mathcal{R}ic^{SB,\mathbb R}$ to $T_{\mathbb{C}}M$. The $(1,1)$-component will be regarded as the associated $(1,1)$-form
\[
\mathcal{R}ic^{SB,\mathbb C}_{(1,1)}
=
\sqrt{-1}\mathcal{R}ic^{SB,\mathbb C}\big(\frac{\partial}{\partial z^i},\frac{\partial}{\partial\bar{z}^j}\big)dz^i\wedge d\bar z^j.
\]
In what follows, we set
\[
\mathrm{Re}\{\mathcal{R}ic^{SB,\mathbb{C}}_{(1,1)}\}=\frac{1}{2}(\mathcal{R}ic^{SB,\mathbb C}_{(1,1)}+\overline{\mathcal{R}ic^{SB,\mathbb C}_{(1,1)}}).
\]
By contrast, the $(2,0)$- and $(0,2)$-components are regarded as complex bilinear tensors:
\[
\mathcal{R}ic^{SB,\mathbb C}_{(2,0)}
=
\mathcal{R}ic^{SB,\mathbb C}\big(\frac{\partial}{\partial z^i},\frac{\partial}{\partial z^j}\big)dz^i\otimes dz^j,
\]
\[
\mathcal{R}ic^{SB,\mathbb C}_{(0,2)}
=
\mathcal{R}ic^{SB,\mathbb C}\big(\frac{\partial}{\partial\bar{z}^i},\frac{\partial}{\partial\bar{z}^j}\big)d\bar z^i\otimes d\bar z^j.
\]

Independently of the $(1,1)$-component of the complexified real Ricci tensor $\mathcal{R}ic^{SB,\mathbb{C}}_{(1,1)}$ that is defined above, the curvature tensor of \(\nabla^{SB}\) on \((T^{1,0}M,h)\) admits four natural Ricci contractions of type \((1,1)\). We denote the first, second, third and fourth Strominger-Bismut-Ricci curvatures by
\[Ric^{SB(1)}=\sqrt{-1}R^{SB(1)}_{i\bar{j}}dz^i\wedge d\bar{z}^{j}\ \ \text{with}\ \ R^{SB(1)}_{i\bar{j}}=h^{k\bar{l}}R^{SB,\mathbb{C}}_{i\bar{j}k\bar{l}},\]
\[Ric^{SB(2)}=\sqrt{-1}R^{SB(2)}_{i\bar{j}}dz^i\wedge d\bar{z}^{j}\ \ \text{with}\ \ R^{SB(2)}_{i\bar{j}}=h^{k\bar{l}}R^{SB,\mathbb{C}}_{k\bar{l}i\bar{j}},\]
\[Ric^{SB(3)}=\sqrt{-1}R^{SB(3)}_{i\bar{j}}dz^i\wedge d\bar{z}^{j}\ \ \text{with}\ \ R^{SB(3)}_{i\bar{j}}=h^{k\bar{l}}R^{SB,\mathbb{C}}_{i\bar{l}k\bar{j}},\]
and
\[Ric^{SB(4)}=\sqrt{-1}R^{SB(4)}_{i\bar{j}}dz^i\wedge d\bar{z}^{j}\ \ \text{with}\ \ R^{SB(4)}_{i\bar{j}}=h^{k\bar{l}}R^{SB,\mathbb{C}}_{k\bar{j}i\bar{l}},\]
respectively, where $R^{SB,\mathbb{C}}_{i\bar{j}k\bar{l}}=R^{SB,\mathbb{C}}(\frac{\partial}{\partial z^i},\frac{\partial}{\partial \bar{z}^j},\frac{\partial}{\partial z^k},\frac{\partial}{\partial \bar{z}^l})$ are the components of the ($\mathbb{C}$-linear) complexified curvature tensor of $^{SB}\nabla$.

In this paper, we collectively refer to the various types of Ricci curvatures associated with the Strominger-Bismut connection $^{SB}\nabla$ on $(M,\omega)$ as \textit{Strominger-Bismut-Ricci curvatures}.

Set $$S_{SB(1)}:=\mathrm{tr}_\omega Ric^{SB(1)}=\mathrm{tr}_\omega Ric^{SB(2)},$$
and
$$S_{SB(2)}:=\mathrm{tr}_\omega Ric^{SB(3)}=\mathrm{tr}_\omega Ric^{SB(4)}$$
be the first and second scalar curvature of $^{SB}\nabla$, respectively.

The main theorems of this paper are below.

\begin{theorem}\label{thm1.2}
Let $(M,\omega)$ be a compact Hermitian surface. If $\mathcal{R}ic^{SB,\mathbb{C}}_{(2,0)}=0$ and
\begin{equation}\label{1.2}
Ric^{SB(2)}+\frac{5}{2}\sqrt{-1}\bar{\partial}^*\omega\wedge\partial^*\omega\leq0,
\end{equation}
then $(M,\omega)$ is a K\"ahler surface.
\end{theorem}

Note that neither $Ric^{SB(3)}$ nor $Ric^{SB(4)}$ is Hermitian symmetric, whereas the sum $Ric^{SB(3)}+Ric^{SB(4)}$ is. It makes sense to define its (semi-)positivity.
\begin{theorem}\label{thm1.3}
Let $(M,\omega)$ be a compact Hermitian surface. Suppose that $\mathcal{R}ic^{SB,\mathbb{C}}_{(2,0)} = 0$,
and that either
\begin{equation}\label{1.3}
Ric^{SB(3)} + Ric^{SB(4)} + 2\sqrt{-1}\bar{\partial}^*\omega \wedge \partial^*\omega \leq 0
\end{equation}
or
\begin{equation}\label{1.4}
\mathrm{Re}\{\mathcal{R}ic^{SB,\mathbb{C}}_{(1,1)}\} + \sqrt{-1}\bar{\partial}^*\omega \wedge \partial^*\omega \leq 0,
\end{equation}
then $(M,\omega)$ is a K\"ahler surface.
\end{theorem}

Streets-Tian \cite{ST10} defined that a Hermitian metric $\omega$ is Hermitian-symplectic if there exists a $(2,0)$-form $\alpha$ such that $d(\alpha+\omega+\bar{\alpha})=0$, and proved that a compact Hermitian surface is Hermitian-symplectic if and only if it is K\"ahler (see another proof in \cite[Theorem 1.2]{LZ09}). A Hermitian-symmetric metric must be pluriclosed, namely, $\partial\bar{\partial}\omega=0$, which in complex dimension two is equivalent to that $\omega$ is Gauduchon.

Every compact complex surface admits a Gauduchon metric (see \cite{Gau77}). When $\omega$ in Theorems \ref{thm1.2} and \ref{thm1.3} is assumed to be a Gauduchon metric, the non-positivity conditions on the second and third Strominger-Bismut-Ricci curvatures can be further relaxed, respectively.
\begin{theorem}\label{thm1.4x}
Let $(M,\omega)$ be a compact Hermitian surface, and suppose that $\omega$ is a Gauduchon metric. If $\mathcal{R}ic^{SB,\mathbb{C}}_{(2,0)}=0$, and
\begin{equation}
Ric^{SB(2)}+\frac{1}{2}\sqrt{-1}\bar{\partial}^*\omega\wedge\partial^*\omega\leq0
\end{equation}
then $(M,\omega)$ is a K\"ahler surface.
\end{theorem}

\begin{theorem}\label{thm1.4}
Let $(M,\omega)$ be a compact Hermitian surface with $\omega$ Gauduchon. Suppose that $\mathcal{R}ic^{SB,\mathbb{C}}_{(2,0)} = 0$
and that either
\begin{equation}\label{1.5}
Ric^{SB(3)} + Ric^{SB(4)} + \sqrt{-1}\,\bar{\partial}^*\omega \wedge \partial^*\omega \leq 0
\end{equation}
or
\begin{equation}\label{1.6}
\mathrm{Re}\{\mathcal{R}ic^{SB,\mathbb{C}}_{(1,1)}\} + \frac{1}{2}\sqrt{-1}\,\bar{\partial}^*\omega \wedge \partial^*\omega \leq 0,
\end{equation}
then $(M,\omega)$ is a K\"ahler surface.
\end{theorem}

\begin{remark}
The assumption \(Ric^{SB,C}_{(2,0)}=0\) plays an essential role in
the arguments below. In particular, statements involving merely
parallel Strominger-Bismut torsion and semi-definiteness of the
Bismut Ricci contractions are not sufficient to force the underlying
surface to be K\"ahler. Standard Hopf surfaces endowed with their
Vaisman metrics provide a basic obstruction to such a conclusion.
\end{remark}

This paper is organized as follows. In Section \ref{sec2}, we fix the notation and present some preliminary lemmas. In Section \ref{sec3}, we establish several identities involving the Ricci curvatures and torsion of the Strominger-Bismut connection on compact Hermitian surfaces. Section \ref{sec4} is devoted to deriving Chern number identities for the Strominger-Bismut-Ricci curvatures. In Section \ref{sec5}, we apply these identities to complete the proofs of Theorems \ref{thm1.2} to \ref{thm1.4}. Finally, in Section \ref{sec6}, we prove certain K\"ahlerness theorems under boundedness conditions on the complexification of the real Strominger-Bismut-Ricci curvatures.
\section{Preliminaries}
\label{sec2}
Let $\{z^1,z^2\}$ be the local holomorphic coordinates on the Hermitian surface $M$, and let $\{x^1,x^2,x^3,x^4\}$ be the associated local real coordinates on the underlying Riemannian manifold $(M,g,J)$ with
\[z^1=x^1+\sqrt{-1}x^3,\ \ \text{and}\ \ z^2=x^2+\sqrt{-1}x^4.\]

Let $(T^{1,0}M,h)$ be the Hermitian holomorphic tangent bundle.  The Chern connection ${^C}\nabla$ is the unique affine connection which is compatible with the Hermitian metric and the holomorphic structure. The Chern connection coefficients are given by
\begin{eqnarray}
{^C}\Gamma_{ij}^k=h^{k\bar{l}}\frac{\partial h_{j\bar{l}}}{\partial z^i},\ \  {^C}\Gamma_{\bar{i}j}^k={^C}\Gamma_{\bar{i}j}^{\bar{k}}={^C}\Gamma_{ij}^{\bar{k}}=0,
\end{eqnarray}
and curvature components by
\begin{equation}
\Theta_{i\bar{j}k\bar{l}}=h_{p\bar{l}}\Theta_{i\bar{j}k}^p=-h_{p\bar{l}}\frac{\partial{^C}\Gamma_{ik}^p}{\partial\bar{z}^j}=-\frac{\partial^2h_{k\bar{l}}}{\partial z^i\partial\bar{z}^j}+h^{p\bar{q}}\frac{\partial h_{p\bar{l}}}{\partial\bar{z}^j}\frac{\partial h_{k\bar{q}}}{\partial z^i}.
\end{equation}
The (first) Chern-Ricci curvature
\begin{equation}
\Theta^{(1)}=\sqrt{-1}\Theta^{(1)}_{i\bar{j}}dz^i\wedge d\bar{z}^j
\end{equation}
represents the first Bott-Chern class $c^{BC}_1(M)$ of $M$, where
\begin{equation}
\Theta^{(1)}_{i\bar{j}}=h^{k\bar{l}}\Theta_{i\bar{j}k\bar{l}}=-\frac{\partial^2\log\det(h_{k\bar{l}})}{\partial z^i\partial\bar{z}^j}.
\end{equation}

The torsion tensor ${^C}T$ of the Chern connection ${^C}\nabla$ on a Hermitian manifold $(M,h)$ is defined by
\begin{equation}
{^C}T_{ij}^k={^C}\Gamma_{ij}^k-{^C}\Gamma_{ji}^k=h^{k\bar{l}}\big(\frac{\partial h_{j\bar{l}}}{\partial z^i}-\frac{\partial h_{i\bar{l}}}{\partial z^j}\big).
\end{equation}
Set
\begin{equation}
T_i=\sum_{k}{^C}T_{ik}^k,\quad\text{and}\quad T_{\bar{i}}=\overline{T_i}.
\end{equation}

The Strominger-Bismut connection $^{SB}\nabla$ is the unique canonical Hermitian connection with totally skew-symmetric torsion (regarded as a 3--form), namely, ${^{SB}}\nabla g=0$, ${^{SB}}\nabla J=0$ and ${^{SB}}T\in\Gamma(M,\wedge^3T^{*}_{\mathbb{R}}M)$ with
\[^{SB}T(X,Y,Z):=g({^{SB}}\nabla_XY-{^{SB}}\nabla_YX-[X,Y],Z)\]
for any $X$, $Y$, $Z\in \Gamma(M,T_{\mathbb{R}}M)$.

The relation between the Levi-Civita connection ${^{LC}}\nabla$ and the Strominger-Bismut connection ${^{SB}}\nabla$ on $(M,g,J)$ is
\begin{equation}
g({^{SB}}\nabla_XY,Z)=g({^{LC}}\nabla_XY,Z)+\frac{1}{2}(d\omega)(JX,JY,JZ)
\end{equation}
for any $X$, $Y$, $Z\in\Gamma(M,T_{\mathbb{R}}M)$. By complexification, it follows that
\begin{equation}
{^{SB}}\Gamma_{\alpha\beta}^{\gamma}={^{LC}}\Gamma_{\alpha\beta}^{\gamma}+\frac{1}{2}{^{SB}}T_{\alpha\beta}^{\gamma}
\end{equation}
with
\begin{equation}
^{LC}\Gamma_{\alpha\beta}^{\gamma}=\frac{1}{2}h^{\gamma\eta}\big(\frac{\partial h_{\alpha\eta}}{\partial z^{\beta}}+\frac{\partial h_{\beta\eta}}{\partial z^{\alpha}}-\frac{\partial h_{\alpha\beta}}{\partial z^{\eta}}\big),
\end{equation}
where $\alpha$, $\beta$, $\gamma$, $\eta\in\{1,2,\bar{1},\bar{2}\}$. Hence, the Strominger-Bismut connection coefficients on $(T^{1,0}M,h)$ are
\begin{equation}\label{2.5}
{^{SB}}\Gamma_{ij}^k=h^{k\bar{l}}\frac{\partial h_{i\bar{l}}}{\partial z^j},
\end{equation}
\begin{equation}\label{2.6}
{^{SB}}\Gamma_{\bar{i}j}^k=h^{k\bar{l}}\big(\frac{\partial h_{j\bar{l}}}{\partial \bar{z}^i}-\frac{\partial h_{j\bar{i}}}{\partial \bar{z}^l}\big),
\end{equation}
and
\begin{equation}\label{2.7}
^{SB}\Gamma_{i\bar{j}}^k={^{SB}}\Gamma_{ij}^{\bar{k}}=0,
\end{equation}
while the torsion tensor ${^{SB}}T$ of the Strominger-Bismut connection ${^{SB}}\nabla$ is
\begin{equation}\label{2.13}
{^{SB}}T_{ij}^k={^{SB}}\Gamma_{ij}^k-{^{SB}}\Gamma_{ji}^k=h^{k\bar{l}}\big(\frac{\partial h_{i\bar{l}}}{\partial z^j}-\frac{\partial h_{j\bar{l}}}{\partial z^i}\big)=-{^{SB}}T_{ji}^k={^{C}}T_{ji}^k
\end{equation}
with
\begin{equation}\label{2.14}
T_i=\sum_{k}{^{SB}}T_{ki}^k=-\sum_{k}{^{SB}}T_{ik}^k.
\end{equation}

By the Bochner formula (see e.g. \cite[Lemma 4.3]{LY12}) that
\[[\bar{\partial}^*,L]=\sqrt{-1}(\partial+[\Lambda,\partial\omega]),\]
it is clear that
\begin{equation}\label{2.22}
\bar{\partial}^*\omega=\sqrt{-1}\Lambda(\partial\omega)=\sqrt{-1}T_idz^i,
\end{equation}
and
\begin{equation}\label{2.22x}
\partial^*\omega=-\sqrt{-1}\Lambda(\bar{\partial}\omega)=-\sqrt{-1}T_{\bar{i}}d\bar{z}^i,
\end{equation}

For any differential forms $\alpha$ and $\beta$ of the same bidegree, we denote by $\langle \alpha,\beta \rangle$ their pointwise inner product and $|\alpha|^2=\langle\alpha,\alpha\rangle$. Define
\[(\alpha,\beta):=\int_M\langle\alpha,\beta\rangle\frac{\omega^2}{2}\quad\text{and}\quad\|\alpha\|^2:=(\alpha,\alpha).\]

To establish our framework, we recall several computational lemmas.
\begin{lemma}[see, e.g., Lemma 3.4 in \cite{LY17}]\label{lem2.2}
Let $(M,h)$ be a Hermitian manifold. For any $p\in M$, there exists holomorphic ''normal coordinates'' $\{z^i\}$ centered at $p$ such that
\begin{equation}\label{2.0}
h_{i\bar{j}}(p)=\delta_{ij},\ \ \frac{\partial h_{i\bar{j}}}{\partial z^k}(p)=-\frac{\partial h_{k\bar{j}}}{\partial z^i}(p),\ \ \text{and}\ \ \frac{\partial h_{i\bar{k}}}{\partial\bar{z}^j}(p)=-\frac{\partial h_{i\bar{j}}}{\partial\bar{z}^k}(p).
\end{equation}
\end{lemma}

As shown in \cite[Lemma 2.5]{Yang2507}, in local holomorphic coordinates, the $(1,1)$-component of the complexification of real Ricci curvature of $^{SB}\nabla$ coincides, in terms of component expressions, with either the third or the fourth Strominger-Bismut-Ricci curvature.
\begin{lemma}[\cite{Yang2507}]\label{lem2.3}
Let $(M,\omega)$ be a Hermitian manifold. For any $X$, $Y\in \Gamma(M,T_{\mathbb{C}}M)$, the complexification of real Ricci curvature of $^{SB}\nabla$ defined in \eqref{1.0} is
\begin{equation}\label{2.8}
\mathcal{R}ic^{SB,\mathbb{C}}(X,Y)=h^{i\bar{l}}R^{SB,\mathbb{C}}\big(\frac{\partial}{\partial z^i},X,Y,\frac{\partial}{\partial\bar{z}^l}\big)+h^{l\bar{i}}R^{SB,\mathbb{C}}\big(\frac{\partial}{\partial \bar{z}^i},X,Y,\frac{\partial}{\partial z^l}\big).
\end{equation}
In particular,
\begin{equation}\label{2.9}
\mathcal{R}_{i\bar{j}}^{SB,\mathbb{C}}=\mathcal{R}ic^{SB,\mathbb{C}}\big(\frac{\partial}{\partial z^i},\frac{\partial}{\partial \bar{z}^j}\big)=h^{l\bar{k}}R^{SB,\mathbb{C}}_{i\bar{k}l\bar{j}}=R^{SB(3)}_{i\bar{j}},
\end{equation}
\begin{equation}\label{2.10}
\mathcal{R}_{\bar{i}j}^{SB,\mathbb{C}}=\mathcal{R}ic^{SB,\mathbb{C}}\big(\frac{\partial}{\partial \bar{z}^i},\frac{\partial}{\partial z^j}\big)=h^{k\bar{l}}R^{SB,\mathbb{C}}_{k\bar{i}j\bar{l}}=R^{SB(4)}_{j\bar{i}},
\end{equation}
and
\begin{equation}\label{2.11}
\mathcal{R}_{ij}^{SB,\mathbb{C}}=\mathcal{R}ic^{SB,\mathbb{C}}\big(\frac{\partial}{\partial z^i},\frac{\partial}{\partial z^j}\big)=h^{k\bar{l}}R^{SB,\mathbb{C}}_{kij\bar{l}},
\end{equation}
\begin{equation}\label{2.12}
\mathcal{R}_{\bar{i}\bar{j}}^{SB,\mathbb{C}}=\mathcal{R}ic^{SB,\mathbb{C}}\big(\frac{\partial}{\partial \bar{z}^i},\frac{\partial}{\partial \bar{z}^j}\big)=h^{k\bar{l}}R^{SB,\mathbb{C}}_{\bar{l}\bar{i}\bar{j}k}.
\end{equation}
\end{lemma}

\begin{remark}
The basic symmetry properties of the curvature tensor of $^{SB}\nabla$ are $R^{SB,\mathbb{C}}(X,Y,Z,W)=-R^{SB,\mathbb{C}}(Y,X,Z,W)=-R^{SB,\mathbb{C}}(X,Y,W,Z)$
for any $X$, $Y$, $Z$, $W\in\Gamma(M,T_{\mathbb{C}}M)$. In general, the first Bianchi identity fails to hold for $R^{SB,\mathbb{C}}$, $R^{SB,\mathbb{C}}(X,Y,Z,W)\neq R^{SB,\mathbb{C}}(Z,W,X,Y)$, and $\mathcal{R}_{ij}^{SB,\mathbb{C}}\neq \mathcal{R}_{ji}^{SB,\mathbb{C}}$. But $R_{i\bar{j}}^{SB(3)}=\overline{R_{j\bar{i}}^{SB(4)}}$, $\mathcal{R}_{i\bar{j}}^{SB,\mathbb{C}}=\overline{\mathcal{R}_{\bar{i}j}^{SB,\mathbb{C}}}$ and $\mathcal{R}_{\bar{i}\bar{j}}^{SB,\mathbb{C}}=\overline{\mathcal{R}_{ij}^{SB,\mathbb{C}}}$.
\end{remark}

The expressions of Strominger-Bismut-Ricci curvatures and scalar curvatures of $^{SB}\nabla$ on a Hermitian surface $(M,\omega)$ follows directly by \cite[Corollary 1.8]{WY25}, \cite[Lemmas 3.2,3.3]{Yang2508} and the fact of
\begin{equation}\label{2.15}
|\partial\omega|^2=|*\partial*\omega|^2=|\bar{\partial}^*\omega|^2.
\end{equation}
\begin{lemma}[\cite{WY25,Yang2508}]\label{lem2.5}
Let $(M,\omega)$ be a Hermitian surface. The Strominger-Bismut-Ricci curvatures are given by
\begin{equation}\label{2.16}
Ric^{SB(1)}=\Theta^{(1)}-(\partial\partial^*\omega+\bar{\partial}\bar{\partial}^*\omega),
\end{equation}
\begin{equation}\label{2.17}
Ric^{SB(2)}=\Theta^{(1)}-(\Lambda\bar{\partial}\bar{\partial}^*\omega+|\bar{\partial}^*\omega|^2)\omega+2\sqrt{-1}\bar{\partial}^*\omega\wedge\partial^*\omega,
\end{equation}
\begin{equation}\label{2.18}
Ric^{SB(3)}=\Theta^{(1)}-\bar{\partial}\bar{\partial}^*\omega+(\Lambda\bar{\partial}\bar{\partial}^*\omega-2|\bar{\partial}^*\omega|^2)\omega+\sqrt{-1}\bar{\partial}^*\omega\wedge\partial^*\omega,
\end{equation}
\begin{equation}\label{2.19}
Ric^{SB(4)}=\Theta^{(1)}-\partial\partial^*\omega+(\Lambda\bar{\partial}\bar{\partial}^*\omega-2|\bar{\partial}^*\omega|^2)\omega+\sqrt{-1}\bar{\partial}^*\omega\wedge\partial^*\omega.
\end{equation}
Let $S_{C(1)}$ be the first Chern scalar curvature. The scalar curvatures of $^{SB}\nabla$ are related by
\begin{equation}\label{2.20}
S_{SB(1)}=S_{C(1)}-2\Lambda\bar{\partial}\bar{\partial}^*\omega,
\end{equation}
\begin{equation}\label{2.21}
S_{SB(2)}=S_{C(1)}+\Lambda\bar{\partial}\bar{\partial}^*\omega-3|\bar{\partial}^*\omega|^2.
\end{equation}
\end{lemma}

\section{Identities on the Strominger-Bismut connection}
\label{sec3}
In this section, we prove several identities related to Ricci curvatures and torsion terms of the Strominger-Bismut connection on a compact Hermitian surface.
\begin{lemma}\label{lem3.1}
Let $(M,\omega)$ be a Hermitian surface. Denote
\[
{^{SB}}T_{ik\bar\ell}=h_{p\bar\ell}{^{SB}}T^p_{ik}.
\]
Then
\begin{equation}\label{3.2}
R^{SB,\mathbb{C}}_{kij\bar\ell}
=
{^{SB}}\nabla_{j}{^{SB}}T_{ik\bar\ell}
+
{^{SB}}\Gamma^{\bar q}_{j\bar\ell}{^{SB}}T_{ik\bar q}
+
{^{SB}}T^p_{kj}{^{SB}}T_{pi\bar\ell}
-
{^{SB}}T^p_{ij}{^{SB}}T_{pk\bar\ell}.
\end{equation}
Consequently, on a Hermitian surface,
\begin{equation}\label{3.3}
R^{SB,\mathbb{C}}_{ij}=\mathcal{R}ic^{SB,\mathbb{C}}_{(2,0)}\big(\frac{\partial}{\partial z^i},\frac{\partial}{\partial z^j}\big)
=
-{^{SB}}\nabla_{j}T_i.
\end{equation}
\end{lemma}

\begin{proof}
By the definition of the curvature tensor and by the coefficient
formulae for the Strominger-Bismut connection, we have
\[
R^{SB,C}_{kij\bar\ell}
=
\frac{\partial}{\partial z^j}{^{SB}}T_{ik\bar\ell}
+{^{SB}}\Gamma^p_{jk}{^{SB}}T_{pi\bar\ell}
-{^{SB}}\Gamma^p_{ji}{^{SB}}T_{pk\bar\ell}
+{^{SB}}T^p_{kj}{^{SB}}T_{pi\bar\ell}
-{^{SB}}T^p_{ij}{^{SB}}T_{pk\bar\ell}.
\]
Since
\[
{^{SB}}\nabla_{\frac{\partial}{\partial z^j}}{^{SB}}T_{ik\bar\ell}
=
\frac{\partial}{\partial z^j}{^{SB}}T_{ik\bar\ell}
-{^{SB}}\Gamma^p_{ji}{^{SB}}T_{pk\bar\ell}
-{^{SB}}\Gamma^p_{jk}{^{SB}}T_{ip\bar\ell}
-{^{SB}}\Gamma^{\bar q}_{j\bar\ell}{^{SB}}T_{ik\bar q},
\]
and \({^{SB}}T_{ip\bar\ell}=-{^{SB}}T_{pi\bar\ell}\), it follows that
\[
\frac{\partial}{\partial z^j}{^{SB}}T_{ik\bar\ell}
+{^{SB}}\Gamma^p_{jk}{^{SB}}T_{pi\bar\ell}
-{^{SB}}\Gamma^p_{ji}{^{SB}}T_{pk\bar\ell}
=
{^{SB}}\nabla_{\frac{\partial}{\partial z^j}}{^{SB}}T_{ik\bar\ell}
+{^{SB}}\Gamma^{\bar q}_{j\bar\ell}{^{SB}}T_{ik\bar q}.
\]
This proves \eqref{3.2}.

Taking the trace with respect to \(h^{k\bar\ell}\), we obtain
\[
R^{SB,C}_{ij}
=
h^{k\bar\ell}{^{SB}}\nabla_{j}{^{SB}}T_{ik\bar\ell}
+h^{k\bar\ell}\Gamma^{\bar q}_{j\bar\ell}{^{SB}}T_{ik\bar q}
+h^{k\bar\ell}{^{SB}}T^p_{kj}{^{SB}}T_{pi\bar\ell}
-h^{k\bar\ell}{^{SB}}T^p_{ij}{^{SB}}T_{pk\bar\ell}.
\]
The first term is
\[
h^{k\bar\ell}{^{SB}}\nabla_{j}{^{SB}}T_{ik\bar\ell}
=
{^{SB}}\nabla_{\frac{\partial}{\partial z^j}}(h^{k\bar\ell}{^{SB}}T_{ik\bar\ell})
=
-{^{SB}}\nabla_{\frac{\partial}{\partial z^j}}T_i.
\]

Since
\[
{}^{SB}\Gamma^{\bar q}_{j\bar\ell}
=
h^{p\bar q}\left(\frac{\partial h_{p\bar\ell}}{\partial z^j}
-\frac{\partial h_{j\bar\ell}}{\partial z^p}\right)
=
-h^{p\bar q}\,{}^{SB}T_{jp\bar\ell},
\]
we have
\[
h^{k\bar\ell}{}^{SB}\Gamma^{\bar q}_{j\bar\ell}
{}^{SB}T_{ik\bar q}
=
-h^{k\bar\ell}h^{p\bar q}
{}^{SB}T_{jp\bar\ell}{}^{SB}T_{ik\bar q}.
\]
Since \(\dim_{\mathbb C}M=2\), the last contraction equals
\[
-h^{k\bar\ell}h^{p\bar q}
{}^{SB}T_{jp\bar\ell}{}^{SB}T_{ik\bar q}
=
-T_iT_j.
\]
Similarly,
\[
h^{k\bar\ell}{}^{SB}T^p_{kj}{}^{SB}T_{pi\bar\ell}=T_iT_j,
\qquad
h^{k\bar\ell}{}^{SB}T^p_{ij}{}^{SB}T_{pk\bar\ell}=0.
\]

Hence the quadratic terms cancel and
\[
R^{SB,C}_{ij}=-{^{SB}}\nabla_{\frac{\partial}{\partial z^j}}T_i.
\]
\end{proof}

\begin{proposition}
On a compact Hermitian surface $(M,\omega)$, we have
\begin{eqnarray}\label{3.21}
&&(\partial\partial^*\omega-\bar{\partial}\bar{\partial}^*\omega,\sqrt{-1}\bar{\partial}^*\omega\wedge\partial^*\omega)\notag\\
&=&(\mathcal{R}ic^{SB,\mathbb{C}}_{(2,0)},\bar{\partial}^*\omega\otimes\bar{\partial}^*\omega)-(\mathcal{R}ic^{SB,\mathbb{C}}_{(0,2)},\partial^*\omega\otimes\partial^*\omega).
\end{eqnarray}
and
\begin{eqnarray}\label{3.20}
&&(\partial\partial^*\omega+\bar{\partial}\bar{\partial}^*\omega,\sqrt{-1}\bar{\partial}^*\omega\wedge\partial^*\omega)\notag\\
&=&-2(\Lambda\bar{\partial}\bar{\partial}^*\omega,|\bar{\partial}^*\omega|^2)+\frac{3}{2}(|\bar{\partial}^*\omega|^4,1)\notag\\
&&+\frac{1}{2}\|\mathcal{R}_{ij}^{SB,\mathbb{C}}+\mathcal{R}_{ji}^{SB,\mathbb{C}}-T_iT_j\|^2-\frac{1}{2}\|\mathcal{R}_{ij}^{SB,\mathbb{C}}+\mathcal{R}_{ji}^{SB,\mathbb{C}}\|^2.
\end{eqnarray}
\end{proposition}

\begin{proof}
It follows from \eqref{2.22} and \eqref{2.22x} that
\begin{eqnarray}\label{3.16x}
&&(\partial\partial^*\omega-\bar{\partial}\bar{\partial}^*\omega,\sqrt{-1}\bar{\partial}^*\omega\wedge\partial^*\omega)\notag\\
&=&-\int_Mh^{i\bar{l}}h^{k\bar{j}}\big(\frac{\partial T_{\bar{j}}}{\partial z^i}-\frac{\partial T_{i}}{\partial \bar{z}^j}\big)T_{\bar{l}}T_k\frac{\omega^2}{2}\notag\\
&=&-\int_Mh^{i\bar{l}}h^{k\bar{j}}({^C}\nabla_{\frac{\partial}{\partial z^i}}T_{\bar{j}}-{^C}\nabla_{\frac{\partial}{\partial \bar{z}^j}}T_{i})T_{\bar{l}}T_k\frac{\omega^2}{2}\notag\\
&=&-\int_M\big(h^{i\bar{l}}{^C}\nabla_{\frac{\partial}{\partial z^i}}(h^{k\bar{j}}T_kT_{\bar{j}})T_{\bar{l}}-h^{k\bar{j}}{^C}\nabla_{\frac{\partial}{\partial \bar{z}^j}}(h^{i\bar{l}}T_{i}T_{\bar{l}})T_k\big)\frac{\omega^2}{2}\notag\\
&&+\int_Mh^{i\bar{l}}h^{k\bar{j}}(T_{\bar{j}}T_{\bar{l}}{^C}\nabla_{\frac{\partial}{\partial z^i}}T_{k}-T_iT_k{^C}\nabla_{\frac{\partial}{\partial \bar{z}^j}}T_{\bar{l}})\frac{\omega^2}{2}\notag\\
&=&(\partial|\bar{\partial}^*\omega|^2,\sqrt{-1}\bar{\partial}^*\omega)+(\bar{\partial}|\bar{\partial}^*\omega|^2,\sqrt{-1}\partial^*\omega)\notag\\
&&+\frac{1}{2}\int_Mh^{i\bar{l}}h^{k\bar{j}}T_{\bar{j}}T_{\bar{l}}({^C}\nabla_{\frac{\partial}{\partial z^i}}T_{k}+{^C}\nabla_{\frac{\partial}{\partial z^k}}T_{i}))\frac{\omega^2}{2}\notag\\
&&-\frac{1}{2}\int_Mh^{i\bar{l}}h^{k\bar{j}}T_iT_k({^C}\nabla_{\frac{\partial}{\partial \bar{z}^j}}T_{\bar{l}}+{^C}\nabla_{\frac{\partial}{\partial \bar{z}^l}}T_{\bar{j}})\frac{\omega^2}{2}\notag\\
&=&(|\bar{\partial}^*\omega|^2,\sqrt{-1}(\partial^*\bar{\partial}^*+\bar{\partial}^*\partial^*)\omega)\notag\\
&&+\frac{1}{2}\int_Mh^{i\bar{l}}h^{k\bar{j}}T_{\bar{j}}T_{\bar{l}}({^{SB}}\nabla_{\frac{\partial}{\partial z^k}}T_{i}+{^{SB}}\nabla_{\frac{\partial}{\partial z^i}}T_{k}))\frac{\omega^2}{2}\notag\\
&&-\frac{1}{2}\int_Mh^{i\bar{l}}h^{k\bar{j}}T_iT_k({^{SB}}\nabla_{\frac{\partial}{\partial \bar{z}^l}}T_{\bar{j}}+{^{SB}}\nabla_{\frac{\partial}{\partial \bar{z}^j}}T_{\bar{l}})\frac{\omega^2}{2}.
\end{eqnarray}

It is well known that
\begin{equation}\label{3.19x}
\partial^*\bar{\partial}^*+\bar{\partial}^*\partial^*=0.
\end{equation}

By \eqref{3.3}, we obtain that
\begin{eqnarray}\label{3.19xx}
&&\int_Mh^{i\bar{l}}h^{k\bar{j}}T_{\bar{j}}T_{\bar{l}}({^{SB}}\nabla_{\frac{\partial}{\partial z^k}}T_{i}+{^{SB}}\nabla_{\frac{\partial}{\partial z^i}}T_{k})\frac{\omega^2}{2}\notag\\
&=&-\int_Mh^{i\bar{l}}h^{k\bar{j}}T_{\bar{j}}T_{\bar{l}}(\mathcal{R}_{ik}^{SB,\mathbb{C}}+\mathcal{R}_{ki}^{SB,\mathbb{C}})\frac{\omega^2}{2}\notag\\
&=&2(\mathcal{R}ic^{SB,\mathbb{C}}_{(2,0)},\bar{\partial}^*\omega\otimes\bar{\partial}^*\omega).
\end{eqnarray}
and
\begin{eqnarray}\label{3.19xxx}
&&\int_Mh^{i\bar{l}}h^{k\bar{j}}T_iT_k({^{SB}}\nabla_{\frac{\partial}{\partial \bar{z}^l}}T_{\bar{j}}+{^{SB}}\nabla_{\frac{\partial}{\partial \bar{z}^j}}T_{\bar{l}})\frac{\omega^2}{2}\notag\\
&=&-\int_Mh^{i\bar{l}}h^{k\bar{j}}T_{i}T_{k}(\mathcal{R}_{\bar{j}\bar{l}}^{SB,\mathbb{C}}+\mathcal{R}_{\bar{l}\bar{j}}^{SB,\mathbb{C}})\frac{\omega^2}{2}\notag\\
&=&2(\mathcal{R}ic^{SB,\mathbb{C}}_{(0,2)},\partial^*\omega\otimes\partial^*\omega).
\end{eqnarray}

\eqref{3.21} follows by applying \eqref{3.19x}, \eqref{3.19xx} and \eqref{3.19xxx} to \eqref{3.16x}.

Similar to the computation in \eqref{3.16x}, we can get
\begin{eqnarray}\label{3.16}
&&(\partial\partial^*\omega+\bar{\partial}\bar{\partial}^*\omega,\sqrt{-1}\bar{\partial}^*\omega\wedge\partial^*\omega)\notag\\
&=&(|\bar{\partial}^*\omega|^2,\sqrt{-1}(\partial^*\bar{\partial}^*-\bar{\partial}^*\partial^*)\omega)\notag\\
&&+\frac{1}{2}\int_Mh^{i\bar{l}}h^{k\bar{j}}T_{\bar{j}}T_{\bar{l}}({^{SB}}\nabla_{\frac{\partial}{\partial z^k}}T_{i}+{^{SB}}\nabla_{\frac{\partial}{\partial z^i}}T_{k}))\frac{\omega^2}{2}\notag\\
&&+\frac{1}{2}\int_Mh^{i\bar{l}}h^{k\bar{j}}T_iT_k({^{SB}}\nabla_{\frac{\partial}{\partial \bar{z}^l}}T_{\bar{j}}+{^{SB}}\nabla_{\frac{\partial}{\partial \bar{z}^j}}T_{\bar{l}})\frac{\omega^2}{2}.
\end{eqnarray}

It is proved in \cite[Lemma 4.5]{Yang2508} that
\begin{equation}\label{3.17}
\Lambda\partial\partial^*\omega=\Lambda\bar{\partial}\bar{\partial}^*\omega=|\bar{\partial}^*\omega|^2-\sqrt{-1}\partial^*\bar{\partial}^*\omega.
\end{equation}
By \eqref{3.19x} and \eqref{3.17}, we have
\begin{eqnarray}\label{3.18}
&&(|\bar{\partial}^*\omega|^2,\sqrt{-1}(\partial^*\bar{\partial}^*-\bar{\partial}^*\partial^*)\omega)\notag\\
&=&2(|\bar{\partial}^*\omega|^2,\sqrt{-1}\partial^*\bar{\partial}^*\omega)\notag\\
&=&-2(\Lambda\bar{\partial}\bar{\partial}^*\omega,|\bar{\partial}^*\omega|^2)+2(|\bar{\partial}^*\omega|^4,1).
\end{eqnarray}

Note that
\begin{eqnarray}\label{3.19}
&&\int_Mh^{i\bar{l}}h^{k\bar{j}}T_{\bar{j}}T_{\bar{l}}({^{SB}}\nabla_{\frac{\partial}{\partial z^k}}T_{i}+{^{SB}}\nabla_{\frac{\partial}{\partial z^i}}T_{k})\frac{\omega^2}{2}\notag\\
&&+\int_Mh^{i\bar{l}}h^{k\bar{j}}T_iT_k({^{SB}}\nabla_{\frac{\partial}{\partial \bar{z}^l}}T_{\bar{j}}+{^{SB}}\nabla_{\frac{\partial}{\partial \bar{z}^j}}T_{\bar{l}})\frac{\omega^2}{2}\notag\\
&=&\|{^{SB}}\nabla_{\frac{\partial}{\partial z^i}}T_{j}+{^{SB}}\nabla_{\frac{\partial}{\partial z^j}}T_{i}+T_iT_j\|^2\notag\\
&&-\|{^{SB}}\nabla_{\frac{\partial}{\partial z^i}}T_{j}+{^{SB}}\nabla_{\frac{\partial}{\partial z^j}}T_{i}\|^2-\|T_iT_j\|^2\notag\\
&=&\|\mathcal{R}_{ij}^{SB,\mathbb{C}}+\mathcal{R}_{ji}^{SB,\mathbb{C}}-T_iT_j\|^2\notag\\
&&-\|\mathcal{R}_{ij}^{SB,\mathbb{C}}+\mathcal{R}_{ji}^{SB,\mathbb{C}}\|^2-(|\bar{\partial}^*\omega|^4,1),
\end{eqnarray}
where we used \eqref{3.3}.

Applying \eqref{3.18} and \eqref{3.19} to \eqref{3.16}, we obtain \eqref{3.20}.
\end{proof}

\begin{lemma}
Let $(M,\omega)$ be a compact Hermitian surface, then we have
\begin{eqnarray}\label{3.29}
&&(Ric^{SB(2)},\frac{1}{2}(\partial^*\partial\omega+\bar{\partial}^*\bar{\partial}\omega))\notag\\
&=&-\|\Lambda\bar{\partial}\bar{\partial}^*\omega\|^2+3(\Lambda\bar{\partial}\bar{\partial}^*\omega,|\bar{\partial}^*\omega|^2)-\frac{3}{2}(|\bar{\partial}^*\omega|^4,1)\notag\\
&&-\frac{1}{2}\|\mathcal{R}_{ij}^{SB,\mathbb{C}}+\mathcal{R}_{ji}^{SB,\mathbb{C}}-T_iT_j\|^2\notag\\
&&+\frac{1}{2}\|\mathcal{R}_{ij}^{SB,\mathbb{C}}+\mathcal{R}_{ji}^{SB,\mathbb{C}}\|^2,
\end{eqnarray}
and
\begin{eqnarray}\label{3.29x}
&&(Ric^{SB(3)},\frac{1}{2}(\partial^*\partial\omega+\bar{\partial}^*\bar{\partial}\omega))\notag\\
&=&(Ric^{SB(4)},\frac{1}{2}(\partial^*\partial\omega+\bar{\partial}^*\bar{\partial}\omega))\notag\\
&=&(\mathrm{Re}\{\mathcal{R}ic^{SB,\mathbb{C}}_{(1,1)}\},\frac{1}{2}(\partial^*\partial\omega+\bar{\partial}^*\bar{\partial}\omega))\notag\\
&=&\frac{1}{2}\|\bar{\partial}\bar{\partial}^*\omega\|^2+\frac{1}{2}\|\Lambda\bar{\partial}\bar{\partial}^*\omega\|^2-\frac{3}{4}(|\bar{\partial}^*\omega|^4,1)\notag\\
&&-\frac{1}{4}\|\mathcal{R}_{ij}^{SB,\mathbb{C}}+\mathcal{R}_{ji}^{SB,\mathbb{C}}-T_iT_j\|^2\notag\\
&&+\frac{1}{4}\|\mathcal{R}_{ij}^{SB,\mathbb{C}}+\mathcal{R}_{ji}^{SB,\mathbb{C}}\|^2.
\end{eqnarray}
\end{lemma}

\begin{proof}
Since $\Lambda\omega=2$ and \cite[(4.7)]{Yang2508} that
\begin{equation}\label{3.30}
\partial^*\partial\omega+\bar{\partial}\bar{\partial}^*\omega=(\Lambda\bar{\partial}\bar{\partial}^*\omega)\omega,
\end{equation}
we have
\begin{equation}\label{3.31}
\Lambda\partial^*\partial\omega=\Lambda\bar{\partial}\bar{\partial}^*\omega.\notag
\end{equation}
Together with \eqref{3.17}, we also have
\begin{equation}\label{3.31x}
\Lambda\partial^*\partial\omega=\Lambda\bar{\partial}\bar{\partial}^*\omega=\Lambda\partial\partial^*\omega=\Lambda\bar{\partial}^*\bar{\partial}\omega.
\end{equation}

Therefore,
\begin{eqnarray*}
&&(Ric^{SB(2)},\partial^*\partial\omega)\notag\\
&=&(\partial\Theta^{(1)},\partial\omega)-(\Lambda\bar{\partial}\bar{\partial}^*\omega+|\bar{\partial}^*\omega|^2,\Lambda\partial^*\partial\omega)+2(\sqrt{-1}\bar{\partial}^*\omega\wedge\partial^*\omega,\partial^*\partial\omega)\notag\\
&=&-\|\Lambda\bar{\partial}\bar{\partial}^*\omega\|^2+(|\bar{\partial}^*\omega|^2,\Lambda\partial\partial^*\omega)-2(\sqrt{-1}\bar{\partial}^*\omega\wedge\partial^*\omega,\bar{\partial}\bar{\partial}^*\omega)\notag\\
&=&-\|\Lambda\bar{\partial}\bar{\partial}^*\omega\|^2+(\Lambda\bar{\partial}\bar{\partial}^*\omega,|\bar{\partial}^*\omega|^2)-2(\partial\partial^*\omega,\sqrt{-1}\bar{\partial}^*\omega\wedge\partial^*\omega),
\end{eqnarray*}
where we used \eqref{2.17} in the first equality, \eqref{3.30} and \eqref{3.31x} in the second.

By taking conjugate and using \eqref{3.20}, we have
\begin{eqnarray*}
&&(Ric^{SB(2)},\frac{1}{2}(\partial^*\partial\omega+\bar{\partial}^*\bar{\partial}\omega))\notag\\
&=&-\|\Lambda\bar{\partial}\bar{\partial}^*\omega\|^2+(\Lambda\bar{\partial}\bar{\partial}^*\omega,|\bar{\partial}^*\omega|^2)\notag\\
&&-(\partial\partial^*\omega+\bar{\partial}\bar{\partial}^*\omega,\sqrt{-1}\bar{\partial}^*\omega\wedge\partial^*\omega)\notag\\
&=&-\|\Lambda\bar{\partial}\bar{\partial}^*\omega\|^2+3(\Lambda\bar{\partial}\bar{\partial}^*\omega,|\bar{\partial}^*\omega|^2)-\frac{3}{2}(|\bar{\partial}^*\omega|^4,1)\notag\\
&&-\frac{1}{2}\|\mathcal{R}_{ij}^{SB,\mathbb{C}}+\mathcal{R}_{ji}^{SB,\mathbb{C}}-T_iT_j\|^2+\frac{1}{2}\|\mathcal{R}_{ij}^{SB,\mathbb{C}}+\mathcal{R}_{ji}^{SB,\mathbb{C}}\|^2.
\end{eqnarray*}
This is \eqref{3.29}.

It follows from \eqref{2.18} and \eqref{2.19} that
\begin{equation}\label{3.31}
Ric^{SB(3)}-Ric^{SB(4)}=\partial\partial^*\omega-\bar{\partial}\bar{\partial}^*\omega.
\end{equation}
Note that \cite[(4.9)]{Yang2508} is
\begin{equation}\label{3.35}
(\bar{\partial}\bar{\partial}^*\omega,\partial^*\partial\omega)=-\|\partial\bar{\partial}^*\omega\|^2.
\end{equation}
Then we get
\begin{eqnarray*}
&&(Ric^{SB(3)}-Ric^{SB(4)},\frac{1}{2}(\partial^*\partial\omega+\bar{\partial}^*\bar{\partial}\omega))\notag\\
&=&(\partial\partial^*\omega-\bar{\partial}\bar{\partial}^*\omega,\frac{1}{2}(\partial^*\partial\omega+\bar{\partial}^*\bar{\partial}\omega))\notag\\
&=&-\frac{1}{2}\|\partial\bar{\partial}^*\omega\|^2+\frac{1}{2}\|\partial\bar{\partial}^*\omega\|^2\notag\\
&=&0.
\end{eqnarray*}
Using Lemma \ref{lem2.3}, we obtain
\begin{eqnarray}\label{3.29xx}
(Ric^{SB(3)},\frac{1}{2}(\partial^*\partial\omega+\bar{\partial}^*\bar{\partial}\omega))&=&(Ric^{SB(4)},\frac{1}{2}(\partial^*\partial\omega+\bar{\partial}^*\bar{\partial}\omega))\notag\\
&=&(\mathrm{Re}\{\mathcal{R}ic^{SB,\mathbb{C}}_{(1,1)}\},\frac{1}{2}(\partial^*\partial\omega+\bar{\partial}^*\bar{\partial}\omega)).
\end{eqnarray}

Calculating directly, we have
\begin{eqnarray}\label{3.32}
&&(Ric^{SB(3)},\partial^*\partial\omega)\notag\\
&=&(\partial\Theta^{(1)},\partial\omega)+\|\partial\bar{\partial}^*\omega\|^2+(\Lambda\bar{\partial}\bar{\partial}^*\omega-2|\bar{\partial}^*\omega|^2,\Lambda\partial^*\partial\omega)\notag\\
&&+(\sqrt{-1}\bar{\partial}^*\omega\wedge\partial^*\omega,\partial^*\partial\omega)\notag\\
&=&\|\partial\bar{\partial}^*\omega\|^2+\|\Lambda\bar{\partial}\bar{\partial}^*\omega\|^2-2(\Lambda\bar{\partial}\bar{\partial}^*\omega,|\bar{\partial}^*\omega|^2)\notag\\
&&+(\bar{\partial}^*\bar{\partial}\omega,\sqrt{-1}\bar{\partial}^*\omega\wedge\partial^*\omega)\notag\\
&=&\|\bar{\partial}\bar{\partial}^*\omega\|^2-(\Lambda\bar{\partial}\bar{\partial}^*\omega,|\bar{\partial}^*\omega|^2)-(\partial\partial^*\omega,\sqrt{-1}\bar{\partial}^*\omega\wedge\partial^*\omega),
\end{eqnarray}
where we used \eqref{2.18} and \eqref{3.35}
in the first equality, \eqref{3.31x} in the second, and \eqref{3.30} and \cite[(4.10)]{Yang2508} that
\begin{equation}\label{4.13}
\|\bar{\partial}\bar{\partial}^*\omega\|^2=\|\Lambda\bar{\partial}\bar{\partial}^*\omega\|^2+\|\partial\bar{\partial}^*\omega\|^2
\end{equation} in the last.

By taking conjugate, we have
\begin{eqnarray}\label{3.33}
&&(Ric^{SB(3)},\bar{\partial}^*\bar{\partial}\omega)\notag\\
&=&(\bar{\partial}\Theta^{(1)},\bar{\partial}\omega)-(\bar{\partial}^*\omega,\bar{\partial}^*\bar{\partial}^*\partial\omega)+(\Lambda\bar{\partial}\bar{\partial}^*\omega-2|\bar{\partial}^*\omega|^2,\Lambda\bar{\partial}^*\bar{\partial}\omega)\notag\\
&&+(\sqrt{-1}\bar{\partial}^*\omega\wedge\partial^*\omega,\bar{\partial}^*\bar{\partial}\omega)\notag\\
&=&\|\Lambda\bar{\partial}\bar{\partial}^*\omega\|^2-(\Lambda\bar{\partial}\bar{\partial}^*\omega,|\bar{\partial}^*\omega|^2)-(\bar{\partial}\bar{\partial}^*\omega,\sqrt{-1}\bar{\partial}^*\omega\wedge\partial^*\omega).
\end{eqnarray}

Combining \eqref{3.32} and \eqref{3.33} and using \eqref{3.20}, we have
\begin{eqnarray}\label{3.34}
(Ric^{SB(3)},\frac{1}{2}(\partial^*\partial\omega+\bar{\partial}^*\bar{\partial}\omega))&=&\frac{1}{2}\|\bar{\partial}\bar{\partial}^*\omega\|^2+\frac{1}{2}\|\Lambda\bar{\partial}\bar{\partial}^*\omega\|^2-\frac{3}{4}(|\bar{\partial}^*\omega|^4,1)\notag\\
&&-\frac{1}{4}\|\mathcal{R}_{ij}^{SB,\mathbb{C}}+\mathcal{R}_{ji}^{SB,\mathbb{C}}-T_iT_j\|^2\notag\\
&&+\frac{1}{4}\|\mathcal{R}_{ij}^{SB,\mathbb{C}}+\mathcal{R}_{ji}^{SB,\mathbb{C}}\|^2.
\end{eqnarray}

\eqref{3.29x} follows by combining \eqref{3.29xx} and \eqref{3.34}.
\end{proof}

The identity given in \cite[Theorem 1.5]{Yang2508} can be reformulated with respect to the Strominger-Bismut-Ricci curvatures.
\begin{lemma}\label{lem3.2}
On a compact Hermitian surface $(M,\omega)$, the following identities hold.
\begin{eqnarray}\label{3.8}
&&\|\bar{\partial}\bar{\partial}^*\omega\|^2+\|\Lambda\bar{\partial}\bar{\partial}^*\omega\|^2\notag\\
&=&2(Ric^{SB(2)},\sqrt{-1}\bar{\partial}^*\omega\wedge\partial^*\omega)+6(\Lambda\bar{\partial}\bar{\partial}^*\omega,|\bar{\partial}^*\omega|^2)\notag\\
&&-4(|\bar{\partial}^*\omega|^4,1)+\frac{1}{2}\|\mathcal{R}_{ij}^{SB,\mathbb{C}}+\mathcal{R}_{ji}^{SB,\mathbb{C}}\|^2,
\end{eqnarray}

\begin{eqnarray}\label{3.9}
&&\|\bar{\partial}\bar{\partial}^*\omega\|^2+\|\Lambda\bar{\partial}\bar{\partial}^*\omega\|^2\notag\\
&=&2(Ric^{SB(3)},\sqrt{-1}\bar{\partial}^*\omega\wedge\partial^*\omega)+\frac{3}{2}(|\bar{\partial}^*\omega|^4,1)\notag\\
&&+\frac{1}{2}\|\mathcal{R}^{SB,\mathbb{C}}_{ij}+\mathcal{R}^{SB,\mathbb{C}}_{ji}-T_iT_j\|^2\notag\\
&&-(\mathcal{R}ic^{SB,\mathbb{C}}_{(2,0)},\bar{\partial}^*\omega\otimes\bar{\partial}^*\omega)+(\mathcal{R}ic^{SB,\mathbb{C}}_{(0,2)},\partial^*\omega\otimes\partial^*\omega),
\end{eqnarray}
and
\begin{eqnarray}\label{3.9x}
&&\|\bar{\partial}\bar{\partial}^*\omega\|^2+\|\Lambda\bar{\partial}\bar{\partial}^*\omega\|^2\notag\\
&=&2(Ric^{SB(4)},\sqrt{-1}\bar{\partial}^*\omega\wedge\partial^*\omega)+\frac{3}{2}(|\bar{\partial}^*\omega|^4,1)\notag\\
&&+\frac{1}{2}\|\mathcal{R}^{SB,\mathbb{C}}_{ij}+\mathcal{R}^{SB,\mathbb{C}}_{ji}-T_iT_j\|^2\notag\\
&&+(\mathcal{R}ic^{SB,\mathbb{C}}_{(2,0)},\bar{\partial}^*\omega\otimes\bar{\partial}^*\omega)-(\mathcal{R}ic^{SB,\mathbb{C}}_{(0,2)},\partial^*\omega\otimes\partial^*\omega).
\end{eqnarray}
\end{lemma}

\begin{proof}
We recall Yang's expression (\cite{Yang2508}) for the complexified Ricci tensor of the Levi-Civita connection. The $(1,1)$-component of the complexified Riemannian Ricci curvature is
\begin{equation}\label{3.11}
\mathfrak{R}ic^{(1,1)}=\Theta^{(1)}-\frac{1}{2}(\partial\partial^*\omega+\bar{\partial}\bar{\partial}^*\omega)+\frac{\sqrt{-1}}{2}\bar{\partial}^*\omega\wedge\partial^*\omega+(\Lambda\bar{\partial}\bar{\partial}^*\omega-|\bar{\partial}^*\omega|^2)\omega,
\end{equation}
and the $(2,0)$-component of the complexified Riemannian Ricci curvature satisfies
\begin{eqnarray}\label{3.12}
\mathfrak{R}_{ij}&=&\mathfrak{R}ic^{(2,0)}\big(\frac{\partial}{\partial z^i},\frac{\partial}{\partial z^j}\big)\notag\\
&=&-\frac{1}{2}\big({^{C}}\nabla_{\frac{\partial}{\partial z^j}}T_{i}+{^{C}}\nabla_{\frac{\partial}{\partial z^i}}T_{j}+T_iT_j\big).
\end{eqnarray}

It is proved in \cite[Theorems 1.5, 3.1]{Yang2508} that
\begin{eqnarray}\label{3.10}
&&\|\bar{\partial}\bar{\partial}^*\omega\|^2+\|\Lambda\bar{\partial}\bar{\partial}^*\omega\|^2\notag\\
&=&2(\mathfrak{R}ic^{(1,1)},\sqrt{-1}\bar{\partial}^*\omega\wedge\partial^*\omega)+2\|\mathfrak{R}ic^{(2,0)}\|^2+\frac{1}{2}(|\bar{\partial}^*\omega|^4,1).
\end{eqnarray}

It is clear that
\[{^{C}}\nabla_{\frac{\partial}{\partial z^j}}T_{i}+{^{C}}\nabla_{\frac{\partial}{\partial z^i}}T_{j}={^{SB}}\nabla_{\frac{\partial}{\partial z^i}}T_{j}+{^{SB}}\nabla_{\frac{\partial}{\partial z^j}}T_{i}.\]
Together with \eqref{3.3} and \eqref{3.12}, we obtain the curvature relation that
\begin{equation}\label{3.13}
\mathfrak{R}_{ij}=\frac{1}{2}(\mathcal{R}_{ij}^{SB,\mathbb{C}}+\mathcal{R}_{ji}^{SB,\mathbb{C}}-T_iT_j),
\end{equation}

By \eqref{2.17} and \eqref{3.11}, we have
\begin{equation}\label{3.14}
\mathfrak{R}ic^{(1,1)}=Ric^{SB(2)}+2(\Lambda\bar{\partial}\bar{\partial}^*\omega)\omega-\frac{3}{2}\sqrt{-1}\bar{\partial}^*\omega\wedge\partial^*\omega-\frac{1}{2}(\partial\partial^*\omega+\bar{\partial}\bar{\partial}^*\omega).
\end{equation}
Together with \eqref{3.20} and \eqref{3.13}, we get
\begin{eqnarray}\label{3.15}
&&(\mathfrak{R}ic^{(1,1)},\sqrt{-1}\bar{\partial}^*\omega\wedge\partial^*\omega)\notag\\
&=&(Ric^{SB(2)},\sqrt{-1}\bar{\partial}^*\omega\wedge\partial^*\omega)+2(\Lambda\bar{\partial}\bar{\partial}^*\omega,|\bar{\partial}^*\omega|^2)\notag\\
&&-\frac{3}{2}(|\bar{\partial}^*\omega|^4,1)-\frac{1}{2}(\partial\partial^*\omega+\bar{\partial}\bar{\partial}^*\omega,\sqrt{-1}\bar{\partial}^*\omega\wedge\partial^*\omega)\notag\\
&=&(Ric^{SB(2)},\sqrt{-1}\bar{\partial}^*\omega\wedge\partial^*\omega)+3(\Lambda\bar{\partial}\bar{\partial}^*\omega,|\bar{\partial}^*\omega|^2)-\frac{9}{4}(|\bar{\partial}^*\omega|^4,1)\notag\\
&&-\|\mathfrak{R}ic^{(2,0)}\|^2+\frac{1}{4}\|\mathcal{R}_{ij}^{SB,\mathbb{C}}+\mathcal{R}_{ji}^{SB,\mathbb{C}}\|^2.
\end{eqnarray}
Applying \eqref{3.15} to \eqref{3.10}, we get \eqref{3.8}.

By \eqref{2.18} and \eqref{3.11}, we have
\begin{equation}\label{3.14}
\mathfrak{R}ic^{(1,1)}=Ric^{SB(3)}+|\bar{\partial}^*\omega|^2\omega+\frac{1}{2}(\bar{\partial}\bar{\partial}^*\omega-\partial\partial^*\omega)-\frac{\sqrt{-1}}{2}\bar{\partial}^*\omega\wedge\partial^*\omega.
\end{equation}
Using \eqref{3.21}, we get
\begin{eqnarray}\label{3.22}
&&(\mathfrak{R}ic^{(1,1)},\sqrt{-1}\bar{\partial}^*\omega\wedge\partial^*\omega)\notag\\
&=&(Ric^{SB(3)},\sqrt{-1}\bar{\partial}^*\omega\wedge\partial^*\omega)+\frac{1}{2}(|\bar{\partial}^*\omega|^4,1)\notag\\
&&-\frac{1}{2}(\mathcal{R}ic^{SB,\mathbb{C}}_{(2,0)},\bar{\partial}^*\omega\otimes\bar{\partial}^*\omega)+\frac{1}{2}(\mathcal{R}ic^{SB,\mathbb{C}}_{(0,2)},\partial^*\omega\otimes\partial^*\omega)
\end{eqnarray}

Applying \eqref{3.13} and \eqref{3.22} to \eqref{3.10}, we get \eqref{3.9}.

\eqref{3.9x} follows from \eqref{3.21}, \eqref{3.31} and \eqref{3.9}.
\end{proof}
\section{Chern number identities}
\label{sec4}
The Chern number identities obtained by Yang \cite{Yang2508} can be reformulated in terms of the Strominger-Bismut-Ricci curvatures.
\begin{lemma}\label{lem4.1}
Let $(M,\omega)$ be a compact Hermitian surface. We have a Chern number identity associated with $Ric^{SB(1)}$ that
\begin{equation}\label{3.23}
4\pi^2c_1^2(M)=\|S_{SB(1)}\|^2-\|Ric^{SB(1)}\|^2+2\|\partial\bar{\partial}^*\omega\|^2,
\end{equation}
\end{lemma}

\begin{proof}
It is shown in \cite[Theorem 3.1]{Yang2508} that the second Chern-Ricci curvature is
\begin{equation}\label{3.26}
\Theta^{(2)}=\Theta^{(1)}-(\partial\partial^*\omega+\bar{\partial}\bar{\partial}^*\omega)+(\Lambda\bar{\partial}\bar{\partial}^*\omega)\omega.
\end{equation}
Combining with \eqref{2.16}, we have
\begin{equation}\label{3.27}
\|\Theta^{(2)}\|^2=\|Ric^{SB(1)}\|^2+2(S_{SB(1)},\Lambda\bar{\partial}\bar{\partial}^*\omega)+2\|\Lambda\bar{\partial}\bar{\partial}^*\omega\|^2.
\end{equation}

Therefore, the Chern number identity given in \cite[Theorem 7.5]{Yang2508} is equivalent to
\begin{eqnarray*}
4\pi^2c_1^2(M)&=&(S_{C(1)}^2,1)-\|\Theta^{(2)}\|^2+2\|\bar{\partial}\bar{\partial}^*\omega\|^2-2(S_{C(1)},\Lambda\bar{\partial}\bar{\partial}^*\omega)\notag\\
&=&(|S_{SB(1)}+2\Lambda\bar{\partial}\bar{\partial}^*\omega|^2,1)-\|Ric^{SB(1)}\|^2\notag\\
&&-2(S_{SB(1)},\Lambda\bar{\partial}\bar{\partial}^*\omega)-2\|\Lambda\bar{\partial}\bar{\partial}^*\omega\|^2+2\|\bar{\partial}\bar{\partial}^*\omega\|^2\notag\\
&&-2(S_{SB(1)}+2\Lambda\bar{\partial}\bar{\partial}^*\omega,\Lambda\bar{\partial}\bar{\partial}^*\omega)\notag\\
&=&\|S_{SB(1)}\|^2-\|Ric^{SB(1)}\|^2+2\|\partial\bar{\partial}^*\omega\|^2,
\end{eqnarray*}
where we used \eqref{2.20} and \eqref{3.27} in the second equality and \eqref{4.13} in the last.
\end{proof}

\begin{lemma}
Let $(M,\omega)$ be a compact Hermitian surface. We have a Chern number identity associated with $Ric^{SB(2)}$ that
\begin{eqnarray}\label{4.1}
&&4\pi^2c_1^2(M)\notag\\
&=&\|S_{SB(1)}\|^2-\|Ric^{SB(2)}\|^2-2(S_{SB(1)},\Lambda\bar{\partial}\bar{\partial}^*\omega+|\bar{\partial}^*\omega|^2)\notag\\
&&-2\|\Lambda\bar{\partial}\bar{\partial}^*\omega\|^2+12(\Lambda\bar{\partial}\bar{\partial}^*\omega,|\bar{\partial}^*\omega|^2)+2\|\partial\bar{\partial}^*\omega\|^2-6(|\bar{\partial}^*\omega|^4,1)\notag\\
&&-4\|\mathcal{R}_{ij}^{SB,\mathbb{C}}+\mathcal{R}_{ji}^{SB,\mathbb{C}}-T_iT_j\|^2+3\|\mathcal{R}_{ij}^{SB,\mathbb{C}}+\mathcal{R}_{ji}^{SB,\mathbb{C}}\|^2.
\end{eqnarray}
\end{lemma}

\begin{proof}
It follows from \eqref{2.17} and \eqref{3.26} that
\[\Theta^{(2)}=Ric^{SB(2)}+A\]
with
\begin{eqnarray}\label{4.2}
A=(2\Lambda\bar{\partial}\bar{\partial}^*\omega+|\bar{\partial}^*\omega|^2)\omega-(\partial\partial^*\omega+\bar{\partial}\bar{\partial}^*\omega)-2\sqrt{-1}\bar{\partial}^*\omega\wedge\partial^*\omega.
\end{eqnarray}
Then we have
\begin{eqnarray}\label{4.3}
\|\Theta^{(2)}\|^2&=&\|Ric^{SB(2)}\|^2+(Ric^{SB(2)},A)+(A,Ric^{SB(2)})+\|A\|^2\notag\\
&=&\|Ric^{SB(2)}\|^2+2(Ric^{SB(2)},A)+\|A\|^2.
\end{eqnarray}

Using \eqref{4.2}, \eqref{3.29} and \eqref{3.8}, we obtain
\begin{eqnarray}\label{4.4}
&&(Ric^{SB(2)},A)\notag\\
&=&(S_{SB(1)},2\Lambda\bar{\partial}\bar{\partial}^*\omega+|\bar{\partial}^*\omega|^2)-(Ric^{SB(2)},\partial^*\partial\omega+\bar{\partial}^*\bar{\partial}\omega)\notag\\
&&-2(Ric^{SB(2)},\sqrt{-1}\bar{\partial}^*\omega\wedge\partial^*\omega)\notag\\
&=&(S_{SB(1)},2\Lambda\bar{\partial}\bar{\partial}^*\omega+|\bar{\partial}^*\omega|^2)+\|\Lambda\bar{\partial}\bar{\partial}^*\omega\|^2-(|\bar{\partial}^*\omega|^4,1)-\|\bar{\partial}\bar{\partial}^*\omega\|^2\notag\\
&&+\|\mathcal{R}_{ij}^{SB,\mathbb{C}}+\mathcal{R}_{ji}^{SB,\mathbb{C}}-T_iT_j\|^2-\frac{1}{2}\|\mathcal{R}_{ij}^{SB,\mathbb{C}}+\mathcal{R}_{ji}^{SB,\mathbb{C}}\|^2.
\end{eqnarray}

By \eqref{4.2}, we have
\begin{eqnarray}\label{4.5}
\|A\|^2&=&2\|2\Lambda\bar{\partial}\bar{\partial}^*\omega+|\bar{\partial}^*\omega|^2\|^2+\|\partial\partial^*\omega+\bar{\partial}\bar{\partial}^*\omega\|^2+4\|\sqrt{-1}\bar{\partial}^*\omega\wedge\partial^*\omega\|^2\notag\\
&&-((2\Lambda\bar{\partial}\bar{\partial}^*\omega+|\bar{\partial}^*\omega|^2)\omega,\partial\partial^*\omega+\bar{\partial}\bar{\partial}^*\omega)\notag\\
&&-(\partial\partial^*\omega+\bar{\partial}\bar{\partial}^*\omega,(2\Lambda\bar{\partial}\bar{\partial}^*\omega+|\bar{\partial}^*\omega|^2)\omega)\notag\\
&&-2((2\Lambda\bar{\partial}\bar{\partial}^*\omega+|\bar{\partial}^*\omega|^2)\omega,\sqrt{-1}\bar{\partial}^*\omega\wedge\partial^*\omega)\notag\\
&&-2(\sqrt{-1}\bar{\partial}^*\omega\wedge\partial^*\omega,(2\Lambda\bar{\partial}\bar{\partial}^*\omega+|\bar{\partial}^*\omega|^2)\omega)\notag\\
&&+2(\partial\partial^*\omega+\bar{\partial}\bar{\partial}^*\omega,\sqrt{-1}\bar{\partial}^*\omega\wedge\partial^*\omega)\notag\\
&&+2(\sqrt{-1}\bar{\partial}^*\omega\wedge\partial^*\omega,\partial\partial^*\omega+\bar{\partial}\bar{\partial}^*\omega)\notag\\
&=&8\|\Lambda\bar{\partial}\bar{\partial}^*\omega\|^2+8(\Lambda\bar{\partial}\bar{\partial}^*\omega,|\bar{\partial}^*\omega|^2)+6(|\bar{\partial}^*\omega|^4,1)\notag\\
&&+2\|\bar{\partial}\bar{\partial}^*\omega\|^2+(\partial\partial^*\omega,\bar{\partial}\bar{\partial}^*\omega)+\overline{(\partial\partial^*\omega,\bar{\partial}\bar{\partial}^*\omega)}\notag\\
&&-2(2\Lambda\bar{\partial}\bar{\partial}^*\omega+|\bar{\partial}^*\omega|^2,\Lambda\partial\partial^*\omega+\Lambda\bar{\partial}\bar{\partial}^*\omega)\notag\\
&&-4(2\Lambda\bar{\partial}\bar{\partial}^*\omega+|\bar{\partial}^*\omega|^2,|\bar{\partial}^*\omega|^2)\notag\\
&&+4(\partial\partial^*\omega+\bar{\partial}\bar{\partial}^*\omega,\sqrt{-1}\bar{\partial}^*\omega\wedge\partial^*\omega).
\end{eqnarray}

It is proved in \cite[(4.12)]{Yang2508} that
\begin{equation}\label{4.6}
(\partial\partial^*\omega,\bar{\partial}\bar{\partial}^*\omega)=\|\Lambda\bar{\partial}\bar{\partial}^*\omega\|^2.
\end{equation}

Applying \eqref{3.20}, \eqref{3.31x} and \eqref{4.6} to \eqref{4.5}, we can get
\begin{eqnarray}\label{4.7}
&&\|A\|^2\notag\\
&=&2\|\bar{\partial}\bar{\partial}^*\omega\|^2+2\|\Lambda\bar{\partial}\bar{\partial}^*\omega\|^2-12(\Lambda\bar{\partial}\bar{\partial}^*\omega,|\bar{\partial}^*\omega|^2)+8(|\bar{\partial}^*\omega|^4,1)\notag\\
&&+2\|\mathcal{R}_{ij}^{SB,\mathbb{C}}+\mathcal{R}_{ji}^{SB,\mathbb{C}}-T_iT_j\|^2-2\|\mathcal{R}_{ij}^{SB,\mathbb{C}}+\mathcal{R}_{ji}^{SB,\mathbb{C}}\|^2.
\end{eqnarray}

Applying \eqref{4.4} and \eqref{4.7} to \eqref{4.3}, we obtain

\begin{eqnarray}\label{4.8}
&&\|\Theta^{(2)}\|^2\notag\\
&=&\|Ric^{SB(2)}\|^2+2(S_{SB(1)},2\Lambda\bar{\partial}\bar{\partial}^*\omega+|\bar{\partial}^*\omega|^2)\notag\\
&&+4\|\Lambda\bar{\partial}\bar{\partial}^*\omega\|^2-12(\Lambda\bar{\partial}\bar{\partial}^*\omega,|\bar{\partial}^*\omega|^2)+6(|\bar{\partial}^*\omega|^4,1)\notag\\
&&+4\|\mathcal{R}_{ij}^{SB,\mathbb{C}}+\mathcal{R}_{ji}^{SB,\mathbb{C}}-T_iT_j\|^2-3\|\mathcal{R}_{ij}^{SB,\mathbb{C}}+\mathcal{R}_{ji}^{SB,\mathbb{C}}\|^2.
\end{eqnarray}

By \eqref{2.20}, \eqref{4.13} and \eqref{4.8}, the Chern number identity given in \cite[Theorem 7.5]{Yang2508} can be reformulated as
\begin{eqnarray*}
&&4\pi^2c_1^2(M)\notag\\
&=&(S_{C(1)}^2,1)-\|\Theta^{(2)}\|^2+2\|\bar{\partial}\bar{\partial}^*\omega\|^2-2(S_{C(1)},\Lambda\bar{\partial}\bar{\partial}^*\omega)\notag\\
&=&(|S_{SB(1)}+2\Lambda\bar{\partial}\bar{\partial}^*\omega|^2,1)-\|Ric^{SB(2)}\|^2-2(S_{SB(1)},2\Lambda\bar{\partial}\bar{\partial}^*\omega+|\bar{\partial}^*\omega|^2)\notag\\
&&-4\|\Lambda\bar{\partial}\bar{\partial}^*\omega\|^2+12(\Lambda\bar{\partial}\bar{\partial}^*\omega,|\bar{\partial}^*\omega|^2)-6(|\bar{\partial}^*\omega|^4,1)\notag\\
&&-4\|\mathcal{R}_{ij}^{SB,\mathbb{C}}+\mathcal{R}_{ji}^{SB,\mathbb{C}}-T_iT_j\|^2+3\|\mathcal{R}_{ij}^{SB,\mathbb{C}}+\mathcal{R}_{ji}^{SB,\mathbb{C}}\|^2\notag\\
&&+2\|\bar{\partial}\bar{\partial}^*\omega\|^2-2(S_{SB(1)}+2\Lambda\bar{\partial}\bar{\partial}^*\omega,\Lambda\bar{\partial}\bar{\partial}^*\omega)\notag\\
&=&\|S_{SB(1)}\|^2-\|Ric^{SB(2)}\|^2-2(S_{SB(1)},\Lambda\bar{\partial}\bar{\partial}^*\omega+|\bar{\partial}^*\omega|^2)\notag\\
&&-2\|\Lambda\bar{\partial}\bar{\partial}^*\omega\|^2+12(\Lambda\bar{\partial}\bar{\partial}^*\omega,|\bar{\partial}^*\omega|^2)+2\|\partial\bar{\partial}^*\omega\|^2-6(|\bar{\partial}^*\omega|^4,1)\notag\\
&&-4\|\mathcal{R}_{ij}^{SB,\mathbb{C}}+\mathcal{R}_{ji}^{SB,\mathbb{C}}-T_iT_j\|^2+3\|\mathcal{R}_{ij}^{SB,\mathbb{C}}+\mathcal{R}_{ji}^{SB,\mathbb{C}}\|^2.
\end{eqnarray*}
This is \eqref{4.1}.
\end{proof}

\begin{lemma}\label{lem4.3}
Let $(M,\omega)$ be a compact Hermitian surface. We have Chern number identity associated with $Ric^{SB(3)}$ that
\begin{eqnarray}\label{4.9}
&&4\pi^2c_1^2(M)\notag\\
&=&\|S_{SB(2)}\|^2-\|Ric^{SB(3)}\|^2+2\|\partial\bar{\partial}^*\omega\|^2+4\|\Lambda\bar{\partial}\bar{\partial}^*\omega\|^2+\frac{5}{2}(|\bar{\partial}^*\omega|^4,1)\notag\\
&&+2(S_{SB(2)},|\bar{\partial}^*\omega|^2)-2(S_{SB(2)},\Lambda\bar{\partial}\bar{\partial}^*\omega)-6(\Lambda\bar{\partial}\bar{\partial}^*\omega,|\bar{\partial}^*\omega|^2)\notag\\
&&-\frac{1}{2}\|\mathcal{R}^{SB,\mathbb{C}}_{ij}+\mathcal{R}^{SB,\mathbb{C}}_{ji}-T_iT_j\|^2
\end{eqnarray}
with
\begin{equation}\label{4.9x}
\|Ric^{SB(3)}\|^2=\|Ric^{SB(4)}\|^2=\|\mathcal{R}ic^{SB,\mathbb{C}}_{(1,1)}\|^2.
\end{equation}
\end{lemma}

\begin{proof}
The Chern number identity on $(M,\omega)$ is
\begin{equation}\label{4.17}
4\pi^2c_1^2(M)=\int_M\Theta^{(1)}\wedge\Theta^{(1)}=\int_M(S_{C(1)}^2-|\Theta^{(1)}|^2)\frac{\omega^2}{2}=\|S_{C(1)}\|^2-\|\Theta^{(1)}\|^2.
\end{equation}

By \eqref{2.21}, we have
\begin{eqnarray}\label{4.18}
\|S_{C(1)}\|^2&=&\|S_{SB(2)}\|^2+\|\Lambda\bar{\partial}\bar{\partial}^*\omega\|^2+9(|\bar{\partial}^*\omega|^4,1)\notag\\
&&-2(S_{SB(2)},\Lambda\bar{\partial}\bar{\partial}^*\omega)+6(S_{SB(2)},|\bar{\partial}^*\omega|^2)\notag\\
&&-6(\Lambda\bar{\partial}\bar{\partial}^*\omega,|\bar{\partial}^*\omega|^2).
\end{eqnarray}

\eqref{2.18} gives that
\[\Theta^{(1)}=Ric^{SB(3)}+B,\]
with
\begin{equation}\label{4.10}
B=\bar{\partial}\bar{\partial}^*\omega-(\Lambda\bar{\partial}\bar{\partial}^*\omega)\omega+2|\bar{\partial}^*\omega|^2\omega-\sqrt{-1}\bar{\partial}^*\omega\wedge\partial^*\omega.
\end{equation}
Therefore,
\begin{eqnarray}\label{4.11}
\|\Theta^{(1)}\|^2&=&\|Ric^{SB(3)}\|^2+(Ric^{SB(3)},B)+(B,Ric^{SB(3)})+\|B\|^2\notag\\
&=&\|Ric^{SB(3)}\|^2+2(Ric^{SB(3)},B)+\|B\|^2.
\end{eqnarray}

Note that
\begin{eqnarray}\label{4.12}
&&(Ric^{SB(3)},B)\notag\\
&=&(Ric^{SB(3)},\bar{\partial}\bar{\partial}^*\omega-(\Lambda\bar{\partial}\bar{\partial}^*\omega)\omega)+2(Ric^{SB(3)},|\bar{\partial}^*\omega|^2\omega)\notag\\
&&-(Ric^{SB(3)},\sqrt{-1}\bar{\partial}^*\omega\wedge\partial^*\omega)\notag\\
&=&-(Ric^{SB(3)},\partial^*\partial\omega)+2(S_{SB(2)},|\bar{\partial}^*\omega|^2)-\frac{1}{2}\|\bar{\partial}\bar{\partial}^*\omega\|^2\notag\\
&&-\frac{1}{2}\|\Lambda\bar{\partial}\bar{\partial}^*\omega\|^2+\frac{3}{4}(|\bar{\partial}^*\omega|^4,1)+\frac{1}{4}\|\mathcal{R}^{SB,\mathbb{C}}_{ij}+\mathcal{R}^{SB,\mathbb{C}}_{ji}-T_iT_j\|^2\notag\\
&&-\frac{1}{2}(\mathcal{R}ic^{SB,\mathbb{C}}_{(2,0)},\bar{\partial}^*\omega\otimes\bar{\partial}^*\omega)+\frac{1}{2}(\mathcal{R}ic^{SB,\mathbb{C}}_{(0,2)},\partial^*\omega\otimes\partial^*\omega)\notag\\
&=&-\frac{3}{2}\|\partial\bar{\partial}^*\omega\|^2-2\|\Lambda\bar{\partial}\bar{\partial}^*\omega\|^2+(\Lambda\bar{\partial}\bar{\partial}^*\omega+2S_{SB(2)},|\bar{\partial}^*\omega|^2)\notag\\
&&+\frac{3}{4}(|\bar{\partial}^*\omega|^4,1)+(\partial\partial^*\omega,\sqrt{-1}\bar{\partial}^*\omega\wedge\partial^*\omega)\notag\\
&&+\frac{1}{4}\|\mathcal{R}^{SB,\mathbb{C}}_{ij}+\mathcal{R}^{SB,\mathbb{C}}_{ji}-T_iT_j\|^2\notag\\
&&-\frac{1}{2}(\mathcal{R}ic^{SB,\mathbb{C}}_{(2,0)},\bar{\partial}^*\omega\otimes\bar{\partial}^*\omega)+\frac{1}{2}(\mathcal{R}ic^{SB,\mathbb{C}}_{(0,2)},\partial^*\omega\otimes\partial^*\omega),
\end{eqnarray}
where we used \eqref{3.30} and \eqref{3.9} in the second equality, and \eqref{3.32}, \eqref{4.13} in the last.

It follows from \eqref{3.21} and \eqref{3.20} that
\begin{eqnarray}\label{4.14}
&&(\partial\partial^*\omega,\sqrt{-1}\bar{\partial}^*\omega\wedge\partial^*\omega)\notag\\
&=&\frac{1}{2}(\mathcal{R}ic^{SB,\mathbb{C}}_{(2,0)},\bar{\partial}^*\omega\otimes\bar{\partial}^*\omega)-\frac{1}{2}(\mathcal{R}ic^{SB,\mathbb{C}}_{(0,2)},\partial^*\omega\otimes\partial^*\omega)\notag\\
&&+\frac{1}{4}\|\mathcal{R}_{ij}^{SB,\mathbb{C}}+\mathcal{R}_{ji}^{SB,\mathbb{C}}-T_iT_j\|^2-\frac{1}{4}\|\mathcal{R}_{ij}^{SB,\mathbb{C}}+\mathcal{R}_{ji}^{SB,\mathbb{C}}\|^2\notag\\
&&-(\Lambda\bar{\partial}\bar{\partial}^*\omega,|\bar{\partial}^*\omega|^2)+\frac{3}{4}(|\bar{\partial}^*\omega|^4,1).
\end{eqnarray}

Applying \eqref{4.14} to \eqref{4.12}, we get
\begin{eqnarray}\label{4.15}
&&(Ric^{SB(3)},B)\notag\\
&=&-\frac{3}{2}\|\partial\bar{\partial}^*\omega\|^2-2\|\Lambda\bar{\partial}\bar{\partial}^*\omega\|^2+2(S_{SB(2)},|\bar{\partial}^*\omega|^2)+\frac{3}{2}(|\bar{\partial}^*\omega|^4,1)\notag\\
&&+\frac{1}{2}\|\mathcal{R}^{SB,\mathbb{C}}_{ij}+\mathcal{R}^{SB,\mathbb{C}}_{ji}-T_iT_j\|^2-\frac{1}{4}\|\mathcal{R}_{ij}^{SB,\mathbb{C}}+\mathcal{R}_{ji}^{SB,\mathbb{C}}\|^2.
\end{eqnarray}

Moreover, \eqref{3.9}, \eqref{4.13} and \eqref{4.10} give that
\begin{eqnarray}\label{4.16}
\|B\|^2&=&\|\bar{\partial}\bar{\partial}^*\omega\|^2+5(|\bar{\partial}^*\omega|^4,1)-2(\Lambda\bar{\partial}\bar{\partial}^*\omega,|\bar{\partial}^*\omega|^2)\notag\\
&&-(\partial\partial^*\omega+\bar{\partial}\bar{\partial}^*\omega,\sqrt{-1}\bar{\partial}^*\omega\wedge\partial^*\omega)\notag\\
&=&\|\partial\bar{\partial}^*\omega\|^2+\|\Lambda\bar{\partial}\bar{\partial}^*\omega\|^2+\frac{7}{2}(|\bar{\partial}^*\omega|^4,1)\notag\\
&&-\frac{1}{2}\|\mathcal{R}_{ij}^{SB,\mathbb{C}}+\mathcal{R}_{ji}^{SB,\mathbb{C}}-T_iT_j\|^2\notag\\
&&+\frac{1}{2}\|\mathcal{R}_{ij}^{SB,\mathbb{C}}+\mathcal{R}_{ji}^{SB,\mathbb{C}}\|^2.
\end{eqnarray}

Applying \eqref{4.15} and \eqref{4.16} to \eqref{4.11}
\begin{eqnarray}\label{4.19}
\|\Theta^{(1)}\|^2&=&\|Ric^{SB(3)}\|^2-2\|\partial\bar{\partial}^*\omega\|^2-3\|\Lambda\bar{\partial}\bar{\partial}^*\omega\|^2+4(S_{SB(2)},|\bar{\partial}^*\omega|^2)\notag\\
&&+\frac{13}{2}(|\bar{\partial}^*\omega|^4,1)+\frac{1}{2}\|\mathcal{R}^{SB,\mathbb{C}}_{ij}+\mathcal{R}^{SB,\mathbb{C}}_{ji}-T_iT_j\|^2.
\end{eqnarray}

We conclude \eqref{4.9} by applying \eqref{4.18} and \eqref{4.19} to \eqref{4.17}.

\eqref{4.9x} follows by Lemma \ref{lem2.3}.
\end{proof}
\section{Proof of main theorems}
\label{sec5}
In this section, we prove Theorems \ref{thm1.2} to \ref{thm1.4} by means of the Ricci curvature and Chern number identities obtained above.
\\\\\textbf{Proof of Theorem \ref{thm1.2}.} Since $\mathcal{R}ic^{SB,\mathbb{C}}_{(2,0)}=0$, we have
\begin{equation}\label{5.3}
\|\mathcal{R}_{ij}^{SB,\mathbb{C}}+\mathcal{R}_{ji}^{SB,\mathbb{C}}\|^2=0.
\end{equation}
Note that $Ric^{SB(2)}+\frac{5}{2}\sqrt{-1}\bar{\partial}^*\omega\wedge\partial^*\omega\leq0$, and $\sqrt{-1}\bar{\partial}^*\omega\wedge\partial^*\omega\geq0$. Applying \eqref{5.3} to \eqref{3.8}, we get
\begin{eqnarray}\label{5.4}
&&\|\bar{\partial}\bar{\partial}^*\omega\|^2+\|\Lambda\bar{\partial}\bar{\partial}^*\omega-3|\bar{\partial}^*\omega|^2\|^2\notag\\
&=&2(Ric^{SB(2)},\sqrt{-1}\bar{\partial}^*\omega\wedge\partial^*\omega)+5(|\bar{\partial}^*\omega|^4,1)\notag\\
&=&2(Ric^{SB(2)}+\frac{5}{2}\sqrt{-1}\bar{\partial}^*\omega\wedge\partial^*\omega,\sqrt{-1}\bar{\partial}^*\omega\wedge\partial^*\omega)\notag\\
&\leq&0.
\end{eqnarray}
It follows that $\bar{\partial}\bar{\partial}^*\omega=0$. Together with \eqref{2.15}, we have
\begin{equation}\label{5.5}
\|\partial\omega\|^2=\|\bar{\partial}^*\omega\|^2=(\bar{\partial}\bar{\partial}^*\omega,\omega)=0.
\end{equation}
Then, $(M,\omega)$ is a K\"ahler surface.$\hfill\Box$
\\\\\textbf{Proof of Theorem \ref{thm1.3}.} By \eqref{2.9} and \eqref{2.10}, we obtain
\[\mathcal{R}ic^{SB,\mathbb{C}}_{(1,1)}=\sqrt{-1}\mathcal{R}^{SB,\mathbb{C}}_{i\bar{j}}dz^i\wedge d\bar{z}^j=Ric^{SB(3)}\]
and
\[\overline{\mathcal{R}ic^{SB,\mathbb{C}}_{(1,1)}}=Ric^{SB(4)}.\]
It follows that
\begin{equation}\label{6.0}
2\mathrm{Re}\{\mathcal{R}ic^{SB,\mathbb{C}}_{(1,1)}\}=\mathcal{R}ic^{SB,\mathbb{C}}_{(1,1)}+\overline{\mathcal{R}ic^{SB,\mathbb{C}}_{(1,1)}}=Ric^{SB(3)}+Ric^{SB(4)}.
\end{equation}

Since $\mathcal{R}ic^{SB,\mathbb{C}}_{(2,0)}=0$, we have
\begin{equation}\label{5.6}
\mathcal{R}ic^{SB,\mathbb{C}}_{(0,2)}=0\quad\text{and}\quad\|\mathcal{R}_{ij}^{SB,\mathbb{C}}+\mathcal{R}_{ji}^{SB,\mathbb{C}}-T_iT_j\|^2=(|\bar{\partial}^*\omega|^4,1).
\end{equation}

Summing up \eqref{3.9} and \eqref{3.9x}, and using \eqref{6.0} and \eqref{5.6}, we obtain
\begin{eqnarray}\label{5.9}
&&\|\bar{\partial}\bar{\partial}^*\omega\|^2+\|\Lambda\bar{\partial}\bar{\partial}^*\omega\|^2\notag\\
&=&(Ric^{SB(3)}+Ric^{SB(4)}+2\sqrt{-1}\bar{\partial}^*\omega\wedge\partial^*\omega,\sqrt{-1}\bar{\partial}^*\omega\wedge\partial^*\omega)\notag\\
&=&2(\mathrm{Re}\{\mathcal{R}ic^{SB,\mathbb{C}}_{(1,1)}\}+\sqrt{-1}\bar{\partial}^*\omega\wedge\partial^*\omega,\sqrt{-1}\bar{\partial}^*\omega\wedge\partial^*\omega).
\end{eqnarray}

Applying \eqref{1.3} or \eqref{1.4} to \eqref{5.9}, we can get $\bar{\partial}\bar{\partial}^*\omega=0$. By \eqref{5.5}, we have $\partial\omega=0$ and hence $(M,\omega)$ is a K\"ahler surface.$\hfill\Box$
\\\\\textbf{Proof of Theorem \ref{thm1.4x}.} Since $\omega$ is Gauduchon, it follows from \eqref{3.17} that
\begin{equation}\label{5.7}
\Lambda\bar{\partial}\bar{\partial}^*\omega=|\bar{\partial}^*\omega|^2.
\end{equation}

Applying \eqref{5.3} and \eqref{5.7} to \eqref{3.8}, we have
\[\|\bar{\partial}\bar{\partial}^*\omega\|^2=2(Ric^{SB(2)}+\frac{1}{2}\sqrt{-1}\bar{\partial}^*\omega\wedge\partial^*\omega,\sqrt{-1}\bar{\partial}^*\omega\wedge\partial^*\omega)\leq0.\]

It follows that $\bar{\partial}\bar{\partial}^*\omega=0$ and then $(M,\omega)$ is a K\"ahler surface by \eqref{5.5}.$\hfill\Box$
\\\\\textbf{Proof of Theorem \ref{thm1.4}.} It follows from \eqref{5.9} and \eqref{5.7} that
\begin{eqnarray}\label{5.10x}
&&\|\bar{\partial}\bar{\partial}^*\omega\|^2\notag\\
&=&(Ric^{SB(3)}+Ric^{SB(4)}+\sqrt{-1}\bar{\partial}^*\omega\wedge\partial^*\omega,\sqrt{-1}\bar{\partial}^*\omega\wedge\partial^*\omega)\notag\\
&=&2(\mathrm{Re}\{\mathcal{R}ic^{SB,\mathbb{C}}_{(1,1)}\}+\frac{1}{2}\sqrt{-1}\bar{\partial}^*\omega\wedge\partial^*\omega,\sqrt{-1}\bar{\partial}^*\omega\wedge\partial^*\omega).
\end{eqnarray}

We get $\bar{\partial}\bar{\partial}^*\omega=0$ by applying \eqref{1.5} or \eqref{1.6} to \eqref{5.10x}. Using \eqref{5.5}, we conclude that $(M,\omega)$ is a K\"ahler surface.$\hfill\Box$

\section{K\"ahlerness theorems under boundedness conditions}\label{sec6}

In this section, we show that a compact Hermitian surface must be K\"ahler if the complexified real Ricci curvature of the Strominger-Bismut connection satisfies appropriate boundedness conditions.

\eqref{2.15} shows that $\bigl(|\bar{\partial}^*\omega|^4,1\bigr)=0$ if and only if $\omega$ is K\"ahler. In particular, when $\omega$ is K\"ahler, we clearly have $\mathcal{R}ic^{SB,\mathbb{C}}_{(2,0)}=0$ and ${^{SB}}T=0$.

If \(\omega\) is not K\"ahler, then \((|\bar\partial^*\omega|^4,1)>0\), and
there exists a positive constant \(a\) such that
\begin{equation}\label{6.1}
\|R^{SB,C}_{ij}+R^{SB,C}_{ji}-T_iT_j\|^2
\leq
a\bigl(|\bar\partial^*\omega|^4,1\bigr).
\end{equation}
If \(\omega\) is K\"ahler, we set \(a=0\).

\begin{theorem}
Let $(M,\omega)$ be a compact Hermitian surface. If
\begin{equation}\label{6.4}
\mathrm{Re}\{\mathcal{R}ic^{SB,\mathbb{C}}_{(1,1)}\}+\frac{a+3}{4}\sqrt{-1}\bar{\partial}^*\omega\wedge\partial^*\omega\leq0,
\end{equation}
then $(M,\omega)$ is a K\"ahler surface.
\end{theorem}

\begin{proof}
It follows from \eqref{3.9}, \eqref{3.9x}, \eqref{6.0} and \eqref{6.1} that
\begin{eqnarray}\label{6.2}
&&\|\bar{\partial}\bar{\partial}^*\omega\|^2+\|\Lambda\bar{\partial}\bar{\partial}^*\omega\|^2\notag\\
&\leq&(2\mathrm{Re}\{\mathcal{R}ic^{SB,\mathbb{C}}_{(1,1)}\},\sqrt{-1}\bar{\partial}^*\omega\wedge\partial^*\omega)+\frac{a+3}{2}(|\bar{\partial}^*\omega|^4,1)\notag\\
&=&2(\mathrm{Re}\{\mathcal{R}ic^{SB,\mathbb{C}}_{(1,1)}\}+\frac{a+3}{4}\sqrt{-1}\bar{\partial}^*\omega\wedge\partial^*\omega,\sqrt{-1}\bar{\partial}^*\omega\wedge\partial^*\omega).\notag\\
\end{eqnarray}

Applying \eqref{6.4} to \eqref{6.2}, we get $\|\bar{\partial}\bar{\partial}^*\omega\|^2=0$ and then $\partial\omega=0$ by \eqref{5.5}. It follows that $(M,\omega)$ is a K\"ahler surface.
\end{proof}

When $\omega$ is Gauduchon, the non-positivity assumption on $\mathcal{R}ic^{SB,\mathbb{C}}_{(1,1)}$ can be significantly relaxed.
\begin{theorem}
Let $(M,\omega)$ be a compact Hermitian surface. If $\omega$ is Gauduchon and
\begin{equation}\label{6.5}
\mathrm{Re}\{\mathcal{R}ic^{SB,\mathbb{C}}_{(1,1)}\}+\frac{a+1}{4}\sqrt{-1}\bar{\partial}^*\omega\wedge\partial^*\omega\leq0,
\end{equation}
then $(M,\omega)$ is a K\"ahler surface.
\end{theorem}

\begin{proof}
Applying \eqref{5.7} to \eqref{6.5}, we get
\begin{eqnarray}
&&\|\bar{\partial}\bar{\partial}^*\omega\|^2\notag\\
&\leq&2(\mathrm{Re}\{\mathcal{R}ic^{SB,\mathbb{C}}_{(1,1)}\}+\frac{a+1}{4}\sqrt{-1}\bar{\partial}^*\omega\wedge\partial^*\omega,\sqrt{-1}\bar{\partial}^*\omega\wedge\partial^*\omega)\notag\\
\end{eqnarray}

Using the condition \eqref{6.5}, we conclude $(M,\omega)$ is a K\"ahler surface as above.
\end{proof}

\section*{Acknowledgements}
The work was supported by Scientific Research Foundation of Chongqing University of Technology (Grant No. 2026ZDZ012). The author is deeply grateful to Prof. Kefeng Liu for his invaluable encouragement and constant guidance, and to Prof. Xiaokui Yang, Prof. Quanting Zhao and Prof. Fangyang Zheng for their insightful discussions and generous support. The author would also like to thank anonymous readers for bringing to his attention an inaccuracy in Lemma \ref{lem3.1} of an earlier version of the paper.

\end{document}